\documentclass[12pt]{article} 

\usepackage{umlaut}

\usepackage[leqno]{amsmath}
\usepackage{amssymb}
\usepackage{amsthm}
\usepackage{pstricks}
\usepackage{enumerate}
\usepackage{color}
\newtheorem{theorem}{Theorem}
\newtheorem{lemma}{Lemma}
\newtheorem{prop}{Proposition}
\newcommand{\mo}{\:\operatorname{mod}\:}
\newcommand{\ph}{\varphi}
\newcommand{\hb}{\frac{1}{2}}
\newcommand{\h}{1/2}

\newcommand{\loglog}{\log \log}
\newcommand{\plm}{\pm}
\newcommand{\zeile}{\vspace{\baselineskip}}

\newcommand{\N}{\mathbb{N}}

\newcommand{\Z}{\mathbb{Z}}
\newcommand{\R}{\mathbb{R}}
\newcommand{\C}{\mathbb{C}}

\allowdisplaybreaks

\begin{document}

\begin{center}
\begin{large}
\textsc{On the Ternary Goldbach Problem with Primes in independent Arithmetic
  Progressions}
\end{large}
\end{center}
\begin{center}
Karin Halupczok
\end{center}

\begin{abstract}
  We show that for every fixed $A>0$ and $\theta>0$ there is a
  $\vartheta=\vartheta(A,\theta)>0$ with the following property.
  Let $n$ be odd and sufficiently large, and let
  $Q_{1}=Q_{2}:=n^{\h}(\log n)^{-\vartheta}$ and $Q_{3}:=(\log
  n)^{\theta}$. 
  Then for all $q_{3}\leq Q_{3}$, all reduced residues $a_{3}$ mod
  $q_{3}$, almost all $q_{2}\leq Q_{2}$, all admissible residues
  $a_{2}$ mod $q_{2}$, almost all $q_{1}\leq Q_{1}$ and all admissible 
  residues $a_{1}$ mod $q_{1}$, there exists a representation
  $n=p_{1}+p_{2}+p_{3}$ with primes $p_{i}\equiv a_{i}\ (q_{i})$,
  $i=1,2,3$.
%
\end{abstract}

%
\subsection{Introduction and results}

\subsubsection{Preliminaries}

Let $n$ be a sufficiently large integer, and for every $i=1,2,3$ let
$a_{i},q_{i}$ be relatively prime integers with $q_{i}\geq 1$ 
and $0\leq a_{i}<q_{i}$.

We consider the ternary Goldbach problem of writing $n$ as
\begin{equation*}
  n=p_{1}+p_{2}+p_{3}
\end{equation*}
with primes $p_{1}$, $p_{2}$ and $p_{3}$ satisfying the three congruences
\begin{equation*}
   p_{i}\equiv a_{i}\mo q_{i}, \; i=1,2,3.
\end{equation*}

A necessary condition for solvability is
\begin{equation*}
  n\equiv a_{1}+a_{2}+a_{3} \mo (q_{1},q_{2},q_{3}),
\end{equation*}
where $(q_{1},q_{2},q_{3})$ denotes the greatest common divisor of the $q_{i}$.
Otherwise no such representation of $n$ is possible.

We precise our consideration in the following way. Let
\begin{equation*}
  J_{3}(n):= \sum_{\substack{ m_{1}+m_{2}+m_{3}=n \\ 
                              m_{i}\equiv a_{i}\: (q_{i}),\\i=1,2,3}} 
    \Lambda (m_{1})\, \Lambda (m_{2})\, \Lambda (m_{3}),
\end{equation*}
where $\Lambda$ is von Mangoldt's function. $J_{3}(n)$ goes closely
with the number of representations of $n$ in the way mentioned.

In this paper we prove that the deviation of $J_{3}(n)$
from its expected main term is
uniformly small for large moduli, namely:
\zeile
\begin{theorem}
\label{Th1}
   For every fixed  $A>0$ and $\theta>0$ there is a
   $\vartheta=\vartheta(A,\theta)>0$ such that 
   for all $q_{3}\leq (\log n)^{\theta}$ and $a_{3}$ with
   $(a_{3},q_{3})=1$ we have
   \begin{multline*}  
     \sum_{q_{2}\leq\frac{ n^{\h}}{(\log n)^{\vartheta}}}
   \; \max_{\substack{a_{2}\\(a_{2},q_{2})=1}} \;
  \sum_{q_{1}\leq \frac{ n^{\h}}{(\log n)^{\vartheta}}}\; 
  \max_{\substack{a_{1}\\(a_{1},q_{1})=1}} 
          \left|  J_{3}(n)   - \frac{n^{2}\mathcal{S}_{3}(n)}
  {2\ph(q_{1})\ph(q_{2})\ph(q_{3})} \right|  \\
\ll \frac{n^{2}}{(\log n)^{A}}. 
\end{multline*} 
The $O$-constant depends on the parameters $A$ and $\theta$.
\end{theorem}

Here $\mathcal{S}_{3}(n)$ denotes the singular series for this
special Goldbach problem and depends on $a_{i}$ and $q_{i}$ likewise
$J_{3}(n)$ does.

We set $S_{3}(n)=0$ if $n\not\equiv a_{1}+a_{2}+a_{3}
\mo (q_{1},q_{2},q_{3})$, the case where trivially $J_{3}(n)=0$
occurs. Then a summand $=0$ in the formula of Theorem \ref{Th1} is given,
therefore we can assume in the proof without loss of generality that
$n\equiv a_{1}+a_{2}+a_{3} \mo (q_{1},q_{2},q_{3})$ holds. We refer to
this as "general condition",
under this, $\mathcal{S}_{3}(n)$ is defined and investigated
later in paragraphs \ref{estsingseries} and \ref{singseries}.

\textbf{Definition.}
For any given $q_{1},q_{2},q_{3}$ we call a triplet $a_{1},a_{2},a_{3}$
of residues mod $q_{1},q_{2},q_{3}$ \textit{admissible} for
$q_{1},q_{2},q_{3}$, if  $(a_{i},q_{i})=1$ for $i=1,2,3$, if
$n\equiv a_{1}+a_{2}+a_{3} \mo (q_{1},q_{2},q_{3})$ and if
$\mathcal{S}_{3}(n)>0$.\\
For given $q_{3},a_{3},q_{2},a_{2}$ and $q_{1}$ we call $a_{1}$
\textit{admissible}, if $a_{1},a_{2},a_{3}$ is admissible for
$q_{1},q_{2},q_{3}$.
For given $q_{3},a_{3},q_{2}$ we call $a_{2}$ \textit{admissible},
if there exists an admissible $a_{1}$ for every positive integer 
$q_{1}$.

We prove in paragraph \ref{singseries}
\begin{lemma}
  \label{lemstar}
  If $n$ is odd, then for given $q_{3}$, $a_{3}$ with
  $(a_{3},q_{3})=1$ and $q_{2}$ there exists an admissible $a_{2}$
  (such that for every $q_{1}$ there exists an admissible $a_{1}$).
  For even $n$ and given $q_{1},q_{2},q_{3}$ 
  there exists \emph{no} admissible triplet $a_{1},a_{2},a_{3}$.
\end{lemma}


Theorem \ref{Th1} provides

\begin{theorem}
\label{th2}
Let $A,\theta,\vartheta>0$ as above and $n\in\N$ odd and sufficiently large.
Let $Q_{1},Q_{2}:= n^{\h}(\log n)^{-\vartheta}$, $Q_{3}:= (\log n)^{\theta}$.
Then for all $q_{3}\leq Q_{3}$, all $a_{3}$, almost all $q_{2}\leq
Q_{2}$, all admissible $a_{2}$, almost all $q_{1}\leq Q_{1}$ and
all admissible $a_{1}$ there exists a representation
$n=p_{1}+p_{2}+p_{3}$ with primes $p_{i}\equiv a_{i} \:(q_{i})$,
$i=1,2,3$.
Here the number of exceptions for $q_{2}$ is
$\ll Q_{2}(\log n)^{-A}$ resp.\ for $q_{1}$ is $\ll Q_{1}(\log
n)^{-A}$. 
\end{theorem}

Theorem \ref{th2} as corollary of Theorem \ref{Th1}
is proved in section \ref{secsix}.

Theorem \ref{Th1} is shown by the circle method. It seems that it
also should hold with the larger bound $q_{3}\leq
n^{\h}(\log n)^{-\vartheta}$, which is the case on the major arcs.
It is not possible to achieve this on the minor arcs by the given methods.

\textbf{Notation.} We denote by $\ph$, $\mu$, $\Lambda$ and $\tau$ the
functions of Euler, M\"obius, von Mangoldt and the divisor
function. Other occuring functions are given in their context.  By
$q_{i}\sim Q_{i}$ we abbreviate $Q_{i}<q_{i}\leq 2Q_{i}$. By $p$ and
$p_{i}$ we denote primes. As usual, $e(\alpha):=e^{2\pi i \alpha}$ for
$\alpha\in\R$.

\subsubsection{Proceeding by the circle method}

Let $A>0$ and $\theta>0$. Let $R:=(\log n)^{B}$ with
$B=B(A,\theta):=\max\{A+\eta+3,D(8A+2\theta+74)\}$,
where $\eta>0$ is some absolute constant (see end of paragraph 
\ref{estsingseries}), and $D(8A+2\theta+74)>0$ is some constant
depending just on $A$ and $\theta$, its definition is given in the
proof of Lemma \ref{lem4}.
Further let $\vartheta>\max\{A+4B+16,\theta+A+3\}$, so $\vartheta$ depends also on
$A$ and $\theta$.

We define major arcs $\mathfrak{M}\subseteq\mathbb{R}$ by
\begin{equation*}
  \mathfrak{M}:= \bigcup_{q\leq R} \bigcup_{\substack{0<a<q \\ (a,q)=1}}
  \left]\frac{a}{q}-\frac{R}{qn}, \frac{a}{q}+\frac{R}{qn} \right[
\end{equation*}
and minor arcs by
$\mathfrak{m}:=\left]-\frac{R}{n},1-\frac{R}{n}\right[
\;\backslash\;\mathfrak{M}$. 

For $\alpha \in \mathbb{R}$ and $j=1,2,3$ let
\begin{equation*}
S_{j}(\alpha):= \sum_{ \substack{m\leq n\\m\equiv a_{j}\,(q_{j})}}
                    \Lambda (m)  \; e(\alpha m).
\end{equation*}

From the orthogonal relations for $e(\alpha m)$ it follows that
\begin{equation*}
  J_{3}(n) = \int_{-\frac{R}{n}}^{1-\frac{R}{n}} S_{1}(\alpha) 
                S_{2}(\alpha) S_{3}(\alpha)\,
                 e(-n\alpha)\, d\alpha.
\end{equation*}

Analogously, denote for $m\leq n$
\begin{equation*}
  J_{2}(m) := \sum_{\substack{m_{2}+m_{3}=m \\ m_{2} \equiv a_{2}
                  \,(q_{2}) \\ m_{3}\equiv a_{3}\, (q_{3})}}
                  \Lambda (m_{2}) \Lambda (m_{3})  
         =     \int_{-\frac{R}{n}}^{1-\frac{R}{n}}
                  S_{2}(\alpha) S_{3}(\alpha)\,e(-m\alpha)\, d\alpha.
\end{equation*}

By
\begin{equation*}
  J^{\mathfrak{M}}_{3} (n):= \int_{\mathfrak{M}} S_{1}(\alpha) 
                S_{2}(\alpha) S_{3}(\alpha)
                 \,e(-n\alpha)\, d\alpha
\end{equation*}
and
\begin{equation*}
  J^{\mathfrak{M}}_{2} (m):=\int_{\mathfrak{M}}S_{2}(\alpha) S_{3}(\alpha)
                  \,e(-m\alpha)\, d\alpha
\end{equation*}
denote the values of $J_{3}(n)$ and $J_{2}(m)$ on the major arcs
$\mathfrak{M}$ and by
\begin{equation*}
  J_{3}^{\mathfrak{m}}(n):= J_{3} (n) - J_{3}^{\mathfrak{M}} (n), \quad
  J_{2}^{\mathfrak{m}}(m) := J_{2} (m) - J_{2}^{\mathfrak{M}} (m)
\end{equation*}
the values on the minor arcs $\mathfrak{m}$.

Concerning the major arcs we get
\begin{theorem} 
\label{Th2}
For $Q_{1},Q_{2},Q_{3}\leq n^{\h}/(\log n)^{\vartheta}$ we have
  \begin{equation*}
    \mathcal{E}_{Q_{1},Q_{2},Q_{3}}^{\mathfrak{M}}:=
      \sum_{\substack{q_{i}\sim Q_{i},\\i=1,2,3}}
      \max_{\substack{a_{i},(a_{i},q_{i})=1,\\i=1,2,3}}
      \left|J_{3}^{\mathfrak{M}} (n) -
      \frac{n^{2}\mathcal{S}_{3}(n)}{2\ph(q_{1})\ph(q_{2})\ph(q_{3})}\right|
      \ll \frac{n^{2}}{(\log n)^{A+3}}.
  \end{equation*}
\end{theorem}

We prove Theorem \ref{Th2} in the following section \ref{seczwo}.

In section \ref{secthree} a lemma containing a special form  of
Montgomery's sieve is proven.
Section \ref{secfour} delivers a proof of Theorem
\ref{Th1} using Theorem \ref{Th2} and the lemma
from section \ref{secthree}.
Further used lemmas concerning estimations on the minor arcs are proven
afterwards in section \ref{secfive}.

\subsection{Estimations on the major arcs} 
\label{seczwo}

\subsubsection{Getting the main term and the error term}
\label{parzwo}

We have
\begin{equation*}
  J_{3}^{\mathfrak{M}} (n) = \sum_{q\leq R} \sum_{\substack{0<a<q\\(a,q)=1}} I(a,q),
\end{equation*}
where
\begin{equation*}
  I(a,q):= \int_{-\frac{R}{qn}}^{\frac{R}{qn}}
  S_{1}\!\left(\frac{a}{q}+\alpha\right) S_{2}\!\left(\frac{a}{q}+\alpha\right) 
  S_{3}\!\left(\frac{a}{q}+\alpha\right)\; e\!\left(-n\left(\frac{a}{q}+\alpha
  \right)\right) d\alpha.
\end{equation*}

For $j=1,2,3$ we have for $\alpha\in[-\frac{R}{qn},\frac{R}{qn}]$
\begin{align*}
  S_{j}\left(\frac{a}{q} + \alpha\right) &= \sum_{\substack{m\leq n \\ m\equiv
  a_{j}\,(q_{j})}} \Lambda(m) \; e(\alpha m)\,
  e\!\left(\frac{a}{q}m\right) \\
  &= \sum_{\substack{m\leq n \\ m\equiv a_{j} (q_{j}) \\ (m,q)=1}}
  \Lambda(m)\; e(\alpha m)\,  e\!\left(\frac{a}{q}m\right)
  + \sum_{\substack{m\leq n \\ m\equiv a_{j} (q_{j}) \\ (m,q)>1}}
  \Lambda(m)\; e(\alpha m)\,  e\!\left(\frac{a}{q}m\right)  \\
  &=  \sum_{\substack{1\leq k\leq q \\ (k,q)=1}} \sum_{\substack{m\leq n \\ m\equiv
  a_{j}\,(q_{j}) \\ m\equiv k\,(q)}} \Lambda(m) \; e(\alpha m)\,
  e\!\left(\frac{a}{q}k\right)  + O((\log n)^{2})
\end{align*}
since
\begin{equation*}
   \sum_{\substack{m\leq n \\ m\equiv a_{j} (q_{j}) \\ (m,q)>1}}
  \Lambda(m) = \sum_{\substack{p^{e}\leq n \\ p^{e}\equiv
  a_{j}(q_{j}) \\ p|q }} \log p \leq \sum_{p|q} \log p \cdot
  \frac{\log n}{\log p} \ll \log n \sum_{p|q} 1  
  \ll (\log n)^{2}.
\end{equation*}

So
\begin{equation*}
   S_{j}\left(\frac{a}{q} + \alpha\right) =
   \sum_{\substack{1\leq k\leq q\\(k,q)=1\\k\equiv a_{j}((q_{j},q))}}
   e\!\left(\frac{a}{q}k\right) 
    T_{j,k}(\alpha)   + O((\log n)^{2})
\end{equation*}
with
\begin{equation*}
  T_{j,k}(\alpha):= \sum_{\substack{m\leq n \\ m\equiv a_{j}\,(q_{j}) \\
  m \equiv k\, (q)}} \Lambda(m) \; e(\alpha m) = \sum_{\substack{m\leq n
  \\ m\equiv f_{j,k}\, ([q_{j},q])}} \Lambda(m)\; e(\alpha m).
\end{equation*}
Here $T_{j,k}$ depends on $k$ with $1\leq k\leq q$, $(k,q)=1$, $k\equiv a_{j}\,((q_{j},q))$.
For such a $k$ there exists an integer $f_{j,k}$
such that the congruence $m\equiv f_{j,k}\,([q_{j},q])$ is equivalent to the system $m\equiv
a_{j}\, (q_{j})$, $m\equiv k \,(q)$, so the last step follows.

Now for positive integers $x$ and $h\leq x$ let 
\begin{equation*}
  \Delta(x,h) := \max_{y\leq x} \max_{(l,h)=1} \Biggl|  
       \sum_{\substack{m\leq y \\ m\equiv l\: (h)}} \Lambda(m)  - \frac{y}{\ph(h)}
     \Biggr|.
\end{equation*}
This expression is $\geq 1$ for $h\leq x$. (Take $y=\ph(h)$ and $l=1$).

Note that by the Theorem of Bombieri and Vinogradov
(see, for example, Br\"udern \cite{c4}) we have
\begin{equation*}
   \sum_{h\leq U} \Delta(x,h) \ll \frac{x}{(\log x)^{D}} + U\sqrt{x}(\log (Ux))^{6}
\end{equation*}
for any fixed $D\geq 1$. This yields that if $U\leq x^{\h}/(\log x)^{D+6}$, then
\begin{equation*}
  \sum_{h\leq U} \Delta(x,h) \ll \frac{x}{(\log x)^{D}}.
\end{equation*}

Now we compute $T_{j,k}(\alpha)$ by partial summation and by introducing
$\Delta$. We get
\begin{align*}
  T_{j,k}(\alpha) &=\sum_{\substack{m\leq n\\ m\equiv f_{j,k}([q_{j},q])
  }} \Lambda(m) \;e(\alpha m) \\
    &= -\int_{0}^{n} \biggl( \sum_{\substack{m\leq y\\ m\equiv
    f_{j,k}([q_{j},q])}} \Lambda(m)\biggr) \frac{d}{dy} (e(\alpha y)) dy
      +  \biggl( \sum_{\substack{m\leq n\\ m\equiv
    f_{j,k}([q_{j},q])}} \Lambda(m)\biggr) e(\alpha n)   \\
   &= -\int_{0}^{n} \biggl(\frac{y}{\ph([q_{j},q])}
                         + O(\Delta(n,[q_{j},q]))\biggr) \frac{d}{dy}
      e(\alpha y) dy  \\
   &\qquad + \biggl( \frac{n}{\ph([q_{j},q])} +
       O(\Delta(n,[q_{j},q]) ) \biggr)e(\alpha n) \\  
 &= \frac{1}{\ph([q_{j},q])} \left( -\int_{0}^{n} y\biggl(\frac{d}{dy}
  e(\alpha y)\biggr) dy + n e(\alpha n) \right) 
   + O\biggl( (1+|\alpha|n) \Delta(n,[q_{j},q])\biggr) \\
 &= \frac{1}{\ph([q_{j},q])} \int_{0}^{n} e(\alpha y) dy +
  O\left(\frac{R}{q}\Delta(n,[q_{j},q])\right), 
\end{align*}
since $|\alpha|\leq \frac{R}{qn}$ and $1\leq \frac{R}{q}$.

This yields, using
\begin{equation*}
  \int_{0}^{n} e(\alpha y) dy = M(\alpha) + O(1),\quad 
  M(\alpha) := \sum_{m=1}^{n} e(\alpha m),
\end{equation*}
the expression
\begin{equation*}
  T_{j,k}(\alpha) = \frac{M(\alpha)}{\ph([q_{j},q])} 
     +O\biggl(\frac{R}{q}\Delta(n,[q_{j},q]) \biggr).
\end{equation*}

We use this term for $T_{j,k}(\alpha)$ to compute $S_{j}(\frac{a}{q}+\alpha)$ as
\begin{align*}
   S_{j}\left(\frac{a}{q}+\alpha\right)  &=
  \sum_{\substack{1\leq k\leq q\\(k,q)=1\\k\equiv a_{j} \:((q_{j},q))}}
   e\!\left(\frac{a}{q}k\right) \left( \frac{M(\alpha)}{\ph([q_{j},q])} 
     +O\left(\frac{R}{q}\Delta(n,[q_{j},q]) \right)
     \right) + O((\log n)^{2}) \\
   &= \frac{c_{j}(a,q)}{\ph([q_{j},q])} M(\alpha) 
    +  O\biggl(\frac{R}{q}\Delta(n,[q_{j},q])\biggr)+O((\log n)^{2})\\
   &= \frac{c_{j}(a,q)}{\ph([q_{j},q])} M(\alpha) 
    +  O\biggl(\frac{R}{q}(\log n)^{2} \Delta(n,[q_{j},q])\biggr)
\end{align*}
since $(\log n)^{2}\geq 1$ and $\frac{R}{q}\Delta(n,[q_{j},q])\geq 1$,
with Ramanujan sums
\begin{equation*}
  c_{j}(a,q):= \sum_{\substack{1\leq k\leq q \\ (k,q)=1 \\ k\equiv
   a_{j} \:((q,q_{j}))}} e\left(\frac{a}{q}k\right) \text{ for } j=1,2,3.
\end{equation*}
We used here that $|c_{j}(a,q)|=1$ or $c_{j}(a,q)=0$, see 
paragraph \ref{estsingseries}.

This provides
\begin{align*}
  I(a,q) &= \int_{-\frac{R}{qn}}^{\frac{R}{qn}}   S_{1} \left(\frac{a}{q}+\alpha \right) S_{2}
    \left(\frac{a}{q}+\alpha\right) S_{3}\left(\frac{a}{q}+\alpha\right)
     e\!\left( -n\left(\frac{a}{q} + \alpha \right) \right) d\alpha \\
    &=H_{a,q}(n)+ \mathcal{O}_{1}+ \mathcal{O}_{2} + \mathcal{O}_{3}
\end{align*}
with
\begin{align*}
  H_{a,q}(n)&:=
    \frac{(c_{1}c_{2}c_{3})(a,q)}{\ph([q_{1},q])\ph([q_{2},q])\ph([q_{3},q])}
    e\left( -n\frac{a}{q}\right) \int_{-\frac{R}{qn}}^{\frac{R}{qn}}
   M^{3}(\alpha) e(-n\alpha)d\alpha,  \\
    \mathcal{O}_{1}&:=  \sum_{j,k,l} \frac{1}{\ph([q_{j},q])\ph([q_{k},q])}
    \int_{-\frac{R}{qn}}^{\frac{R}{qn}} |M^{2}(\alpha)|
      d\alpha  \cdot O\left(  \frac{R}{q} (\log n)^{2}
    \Delta(n,[q_{l},q]) \right),\\ 
   \mathcal{O}_{2}&:=  \sum_{j,k,l} \frac{1}{\ph([q_{j},q])}
    \int_{-\frac{R}{qn}}^{\frac{R}{qn}}  |M(\alpha)|
    d\alpha  \cdot O\left(  \frac{R^{2}}{q^{2}} (\log n)^{4}
    \Delta(n,[q_{k},q])\Delta(n,[q_{l},q]) \right),\\ 
    \mathcal{O}_{3}&:= O\left(  \frac{R^{3}}{q^{3}} (\log n)^{6}
    \Delta(n,[q_{1},q])\Delta(n,[q_{2},q])\Delta(n,[q_{3},q]) \frac{R}{qn}\right).
\end{align*}
Note that we abbreviated
$(c_{1}c_{2}c_{3})(a,q):=c_{1}(a,q)c_{2}(a,q)c_{3}(a,q)$. The sum
$\sum_{j,k,l}$ is over all triplets $(j,k,l)$ of pairwise different $j,k,l\in\{1,2,3\}$.


So we managed to show
\begin{equation*}
  J_{3}^{\mathfrak{M}} (n) = \sum_{q\leq R} \sum_{\substack{a<q \\ (a,q)=1}}
  I(a,q) = \sum_{q\leq R} \sum_{\substack{a<q \\ (a,q)=1}} (H_{a,q}(n) +
  \mathcal{O}_{1} +\mathcal{O}_{2} +\mathcal{O}_{3}).
\end{equation*}

The main term of $J_{3}^{\mathfrak{M}} (n)$ is contained in
\begin{equation*}
  H(n):=\sum_{q\leq R} \sum_{\substack{a<q \\ (a,q)=1}} H_{a,q}(n).
\end{equation*}

We have to show now that for each $i=1,2,3$ the error term
$\mathcal{O}_{i}$ fulfills
\begin{equation*}
  \sum_{q_{1},q_{2},q_{3}} \sum_{q\leq R} \sum_{\substack{a<q\\(a,q)=1}}
  \mathcal{O}_{i} \ll \frac{n^{2}}{(\log n)^{A+3}},
\end{equation*}
then it will follow that
\begin{equation*}
  \sum_{q_{1},q_{2},q_{3}} \max_{a_{1},a_{2},a_{3}}
  \left|J_{3}^{\mathfrak{M}} (n) - H(n)
  \right| \ll \frac{n^{2}}{(\log n)^{A+3}}.
\end{equation*}
The main term $H(n)$ will be considered later.

So we first consider the error term with $\mathcal{O}_{1}$. It is
(since $\ph(q)\gg q/(\loglog q)$)
\begin{align*}
  &\ll \sum_{j,k,l} \sum_{q\leq R} \sum_{q_{j},q_{k}}
  \frac{1}{\ph([q_{j},q])\ph([q_{k},q])} \sum_{\substack{a<q\\(a,q)=1}}
   \frac{R^{2}}{q^{2}}n(\log n)^{2} \sum_{q_{l}}\Delta(n,[q_{l},q]) \\
  &\ll \sum_{j,k,l} \sum_{q_{j}} \frac{\loglog n}{q_{j }} \sum_{q_{k}}
  \frac{\loglog n}{q_{k}} R^{2}n (\log n)^{2} \sum_{q_{l}}
  \sum_{q\leq R} \frac{1}{q} \Delta(n,[q_{l},q]) \\
  &\ll R^{2}n (\log n)^{5} \sum_{j,k,l} \sum_{h_{l}\leq R Q_{l}}
  \omega(h_{l}) \Delta(n,h_{l})
\end{align*}
with
\begin{align*}
  \omega(h_{l}) &:= \sum_{q_{l}} \sum_{\substack{q\leq
  R\\ [q_{l},q]=h_{l}}} \frac{1}{q} = \sum_{d_{l}\leq R} \sum_{q_{l}}
  \sum_{\substack{q\leq R \\ (q_{l},q)=d_{l} \\ q_{l}q = h_{l} d_{l}}}
    \frac{1}{q}\\
    &\ll \sum_{d_{l}\leq R} \sum_{\substack{q\leq R \\ d_{l}|q}} \frac{1}{q}
    \ll \sum_{d_{l}\leq R} \sum_{q\leq R} \frac{1}{qd_{l}} \ll (\log n)^{2},
\end{align*}
so the $\mathcal{O}_{1}$-error term is
\begin{align*}
  &\ll R^{2} n (\log n)^{7} \sum_{j,k,l} \sum_{h\leq RQ_{l}} \Delta(n,h) \ll R^{2}n
  (\log n)^{7} \frac{n}{(\log n)^{D}} \\  &\ll \frac{n^{2}}{(\log
  n)^{D-2B-7}}  \ll \frac{n^{2}}{(\log n)^{A+3}},
\end{align*}
for some $D\geq A+2B+10$ and $D+6\leq\vartheta-B$, so this holds
if $\vartheta\geq A+3B+16$, which is the case. We used
the Theorem of Bombieri and Vinogradov
with $Q_{l}\ll n^{\h}(\log n)^{-\vartheta}$ for $\vartheta>0$. So we
are done for $\mathcal{O}_{1}$.

We consider now the error term with $\mathcal{O}_{2}$. It is
\begin{align*}
  &\ll \sum_{j,k,l}\sum_{q\leq R} \sum_{q_{j}} \frac{1}{\ph([q_{j},q])}
   \sum_{\substack{a<q\\(a,q)=1}} 
   \frac{R^{3}}{q^{3}}(\log n)^{4} \sum_{q_{k},q_{l}}
   \Delta(n,[q_{k},q])\Delta(n,[q_{l},q]) \\
  &\ll \sum_{j,k,l} \sum_{q_{j}} \frac{\loglog n}{q_{j }} R^{3}
   (\log n)^{4} \sum_{q_{k},q_{l}}
  \sum_{q\leq R} \frac{1}{q^{2}} \Delta(n,[q_{k},q])\Delta(n,[q_{l},q]) \\
  &\ll R^{3} (\log n)^{6} \sum_{j,k,l} \sum_{h_{k}\leq R Q_{k}} 
   \sum_{h_{l}\leq R Q_{l}} \omega(h_{k},h_{l}) \Delta(n,h_{k})\Delta(n,h_{l})
\end{align*}
with
\begin{align*}
  \omega(h_{k},h_{l}) &:= \sum_{q_{k},q_{l}} \sum_{\substack{q\leq
  R\\ [q_{k},q]=h_{k} \\ [q_{l},q]=h_{l}}} \frac{1}{q^{2}}
    = \sum_{d_{k},d_{l}\leq R} \sum_{q_{k},q_{l}}
  \sum_{\substack{q\leq R \\ (q_{k},q)=d_{k},(q_{l},q)=d_{l} \\
  q_{k}q=h_{k}d_{k} ,q_{l}q = h_{l} d_{l}}} \frac{1}{q^{2}}\\
    &\ll \sum_{d_{k},d_{l}\leq R} \sum_{\substack{q\leq R \\
  [d_{k},d_{l}]|q}} \frac{1}{q^{2}} 
    \leq \sum_{d_{k},d_{l}\leq R} \sum_{q\leq R}
  \frac{1}{q^{2}[d_{k},d_{l}]^{2}} \\
   &= \sum_{d_{k},d_{l}\leq R} \sum_{q\leq R}
  \frac{(d_{k},d_{l})^{2}}{q^{2}d_{k}^{2}d_{l}^{2}} 
   \leq R^{2} \sum_{q} \frac{1}{q^{2}} \sum_{d_{k}}\frac{1}{d_{k}^{2}}
  \sum_{d_{l}}\frac{1}{d_{l}^{2}} \ll R^{2},
\end{align*}
so the $\mathcal{O}_{2}$-error term is
\begin{align*}
  &\ll R^{5} (\log n)^{6}  \sum_{h_{k}\leq RQ_{k}} \Delta(n,h_{k})
  \sum_{h_{l}\leq RQ_{l}} \Delta(n,h_{l}) \\ &\ll R^{5} (\log n)^{6} 
  \cdot\left(\frac{n}{(\log n)^{D}}\right)^{2}
  = \frac{n^{2}}{(\log n)^{2D-5B-6}} \ll
  \frac{n^{2}}{(\log n)^{A+3}},
\end{align*}
for some $2D\geq A+5B+9$ and $D+6\leq\vartheta-B$, so this holds 
if $\vartheta\geq \frac{1}{2} A+\frac{7}{2}B+11$, which is the case.
We used the Theorem of Bombieri and Vinogradov with
$Q_{k},Q_{l}\ll n^{\h}(\log n)^{-\vartheta}$ for
$\vartheta>0$. So we are done for $\mathcal{O}_{2}$.

Now to the error term with $\mathcal{O}_{3}$, it is
\begin{align*}
  &\ll \sum_{q\leq R} \sum_{\substack{a<q\\ (a,q)=1}}
  \frac{R^{4}}{q^{4}n} (\log n)^{6}  \sum_{q_{1},q_{2},q_{3}}
  \Delta(n,[q_{1},q]) \Delta(n,[q_{2},q]) \Delta(n,[q_{3},q])  \\
  &\ll  \frac{R^{4}}{n} (\log n)^{6} \sum_{\substack{h_{1}\leq RQ_{1}  \\
  h_{2}\leq RQ_{2} \\ h_{3}\leq RQ_{3} }} \omega(h_{1},h_{2},h_{3})
  \Delta(n,h_{1}) \Delta(n,h_{2}) \Delta(n,h_{3}) 
\end{align*}
with
\begin{align*}
  \omega(h_{1},h_{2},h_{3}) &:= \sum_{q_{1},q_{2},q_{3}}
  \sum_{\substack{q\leq R \\ [q_{i},q]=h_{i} }}
  \frac{1}{q^{3}} = \sum_{d_{1},d_{2},d_{3}\leq R} \sum_{q_{1},q_{2},q_{3}}
  \sum_{\substack{q\leq R \\ (q_{i},q)=d_{i} \\ q_{i}q = h_{i}d_{i}}}
  \frac{1}{q^{3}} \ll \sum_{d_{1},d_{2},d_{3}\leq R}
  \sum_{\substack{q\leq R \\ [d_{1},d_{2},d_{3}] | q}} \frac{1}{q^{3}}  \\
   &\ll \sum_{d_{1},d_{2},d_{3}\leq R} \sum_{q\leq R}
  \frac{1}{q^{3}[d_{1},d_{2},d_{3}]^{3}} = \sum_{d_{1},d_{2},d_{3}\leq
  R} \sum_{q\leq R}
  \frac{(d_{1},[d_{2},d_{3}])^{3}(d_{2},d_{3})^{3}}
  {q^{3}d_{1}^{3}d_{2}^{3}d_{3}^{3}} \\ &\ll
  \sum_{d_{1},d_{2}\leq R} \sum_{d_{3}\leq R} \frac{1}{d_{3}^{3}}
  \sum_{q\leq R} \frac{1}{q^{3}} \ll R^{2},
\end{align*}
so the $\mathcal{O}_{3}$-error term is
\begin{align*}
  &\ll \frac{R^{6}}{n} (\log n)^{6} \sum_{h_{1}\leq RQ_{1}}
  \Delta(n,h_{1}) \sum_{h_{2}\leq RQ_{2}} \Delta(n,h_{2}) 
  \sum_{h_{3}\leq RQ_{3}} \Delta(n,h_{3}) \\
  &\ll \frac{R^{6}}{n} (\log n)^{6}
  \frac{n^{3}}{(\log n)^{3D}}
  = \frac{n^{2}}{(\log n)^{3D-6B-6}} \ll
  \frac{n^{2}}{(\log n)^{A}},
\end{align*}
for some $3D\geq A+6B+9$ and $D+6\leq\vartheta-B$, so 
this holds if $\vartheta\geq\frac{1}{3}A+3B+9$, which is the case.
We used the Theorem of Bombieri and Vinogradov with $Q_{1},Q_{2},Q_{3}\ll
n^{\h}(\log n)^{-\vartheta}$ for
$\vartheta>0$. So we are done for $\mathcal{O}_{3}$.

What is now left is the consideration of the main term $H(n)$.
Since
\begin{equation*}
  \int_{-\frac{R}{qn}}^{\frac{R}{qn}} M^{3}(\alpha) e(-n\alpha)d\alpha
   = \frac{n^{2}}{2} + O\left(\frac{q^{2}n^{2}}{R^{2}}\right)
\end{equation*}
(see for example Vaughan \cite{c5}) we have
\begin{equation*}
  H(n)= \sum_{q\leq R} \sum_{\substack{a<q \\ (a,q)=1}}
  \frac{(c_{1}c_{2}c_{3})(a,q) e\bigl(-n\frac{a}{q}\bigr)}
  {\ph([q_{1},q])\ph([q_{2},q])\ph([q_{3},q])} \left(\frac{n^{2}}{2} +
  O\biggl(\frac{q^{2}n^{2}}{R^{2}}\biggr)\right).
\end{equation*}

Now let
\begin{equation*}
  \lambda(q):=
  \frac{\ph(q_{1})\ph(q_{2})\ph(q_{3})}
  {\ph([q_{1},q])\ph([q_{2},q])\ph([q_{3},q])} b(q)  
\end{equation*}
with
\begin{equation*}
  b(q):=\sum_{\substack{a<q \\ (a,q)=1}} (c_{1}c_{2}c_{3})(a,q)\; 
  e\!\left(-n\frac{a}{q}\right)
\end{equation*}
and let
\begin{equation*}
  \mathcal{S}_{3} (n):=  \sum_{q=1}^{\infty} \lambda(q)
\end{equation*}
be the singular series. In the next paragraph we show that
it is absolutely convergent.

Therefore we have
\begin{align*}
  H(n)&= \sum_{q\leq R}
  \frac{\lambda(q)n^{2}}{2\ph(q_{1})\ph(q_{2})\ph(q_{3})} 
   + O\Biggl(\frac{n^{2}}{R^{2}} \sum_{q\leq R}
  \frac{q^{2}|\lambda(q)|}{\ph(q_{1})\ph(q_{2})\ph(q_{3})} \Biggr) \\
    &= \frac{n^{2}}{2\ph(q_{1})\ph(q_{2})\ph(q_{3})}\mathcal{S}_{3}(n) 
   + O(e_{1}) + O(e_{2})
\end{align*}
with
\begin{align*}
  e_{1}&:=\frac{n^{2}}{\ph(q_{1})\ph(q_{2})\ph(q_{3})} \sum_{q>R}
  |\lambda(q)|, \\ 
  e_{2}&:=\frac{n^{2}}{R^{2}\ph(q_{1})\ph(q_{2})\ph(q_{3})}
  \sum_{q\leq R} q^{2}|\lambda(q)| .
\end{align*}
For the two occuring error terms $e_{1}$ and $e_{2}$ we have
to show that
\begin{equation*}
    \sum_{q_{1},q_{2},q_{3}} \max_{a_{1},a_{2},a_{3}} e_{j}
     \ll \frac{n^{2}}{(\log n)^{A+3}},
\end{equation*}
then Theorem \ref{Th2} follows.
This is done in the next paragraph.

\subsubsection{Estimations with the singular series}
\label{estsingseries}

Now we need estimations for the $\lambda$-series.
These show the absolute convergence of $\mathcal{S}_{3}(n)$ and can
also be used to deal with $e_{1}$ and $e_{2}$.

First we state that the Ramanujan sums $c_{j}(a,q)$ for fixed $a_{j},q_{j}$,
$j=1,2,3$, can be computed by
\begin{equation*}
  c_{j}(a,q) = c_{a_{j},q_{j}}(a,q)=\begin{cases}  \mu\big(\frac{q}{d_{j}}\big)
  e\big(\frac{au_{j}a_{j}}{d_{j}}\big), &\text{ if }
  \big(d_{j},\frac{q}{d_{j}}\big)=1 , \\
  0, &\text{ else,} \end{cases}
\end{equation*}
where $d_{j}:=(q_{j},q)$ and $u_{j}$ is the
solution of the congruence $\frac{q}{d_{j}}u_{j} \equiv 1 \:(d_{j})$, with 
$0\leq u_{j}<d_{j}$. 
(For a proof see \cite{c8}). From this result we already used that $|c_{j}(a,q)|=1$ or
$c_{j}(a,q)=0$ in the paragraph before.

We are now going to show that $b$ is multiplicative in $q$. We prove a
proposition about the $c_{j}$ first.

\begin{prop}
\label{Prop1} 
Let $q=\bar{q}\tilde{q}$, $(\bar{q},\tilde{q})=1$, $(a,q)=1$, and
let $a=\tilde{a}\bar{q}+\bar{a}\tilde{q}$ with
$(\tilde{a},\tilde{q})=1$, $(\bar{a},\bar{q})=1$. Then
$c_{j}(a,q)=c_{j}(\tilde{a},\tilde{q})\cdot c_{j}(\bar{a},\bar{q})$
for $j=1,2,3$.
\end{prop}

\textbf{Proof.} 
Let $\tilde{a}_{j}\bar{q}+\bar{a}_{j}\tilde{q}\equiv a_{j}
\:((q_{j},q))$ with $\tilde{a}_{j}$ a residue mod $(q_{j},\tilde{q})$
and $\bar{a}_{j}$ a residue mod $(q_{j},\bar{q})$. Then we have for $j=1,2,3$
\begin{equation*}
  c_{j}(a,q)= \sum_{\substack{m<q\\(m,q)=1\\m\equiv a_{j}((q_{j},q))}}
      e\Bigl(m\frac{a}{q}\Bigr) 
   =
   \sum_{\substack{\tilde{m}<\tilde{q}\\(\tilde{m},\tilde{q})=1\\\tilde{m}\equiv
     \tilde{a}_{j}\:((q_{j},\tilde{q})) }} \:
   \sum_{\substack{\bar{m}<\bar{q}\\(\bar{m},\bar{q})=1\\\bar{m}\equiv
     \bar{a}_{j}\:((q_{j},\bar{q})) }} 
   e\Bigl( \frac{(\tilde{m}\bar{q}+\bar{m}\tilde{q})
   (\tilde{a}\bar{q}+\bar{a}\tilde{q})}{\tilde{q}\bar{q}}  \Bigr)
\end{equation*}
by substituting $m=\tilde{m}\bar{q}+\bar{m}\tilde{q}$ with
$\tilde{m}\equiv \tilde{a}_{j} \:((q_{j},\tilde{q}))$ and
$\bar{m}\equiv \bar{a}_{j} \:((q_{j},\bar{q}))$, and we have
$\tilde{a}_{j} \equiv a_{j}\bar{q}^{-1}\:((q_{j},\tilde{q})) $
and $\bar{a}_{j} \equiv a_{j}\tilde{q}^{-1}\:((q_{j},\bar{q})) $.
Therefore we get
\begin{align*}
  c_{j}(a,q) &=  \sum_{\substack{\tilde{m}<\tilde{q}\\(\tilde{m},\tilde{q})=1\\\tilde{m}\equiv
     a_{j}\bar{q}^{-1}\:((q_{j},\tilde{q})) }}
 e\Bigl(\frac{\tilde{m}\tilde{a}\bar{q}}{\tilde{q}}\Bigr) 
    \sum_{\substack{\bar{m}<\bar{q}\\(\bar{m},\bar{q})=1\\\bar{m}\equiv
     a_{j}\tilde{q}^{-1}\:((q_{j},\bar{q})) }} 
 e\Bigl(\frac{\bar{m}\bar{a}\tilde{q}}{\bar{q}}\Bigr)  \\
  &= \sum_{\substack{\tilde{m}<\tilde{q}\\(\tilde{m},\tilde{q})=1\\\tilde{m}\equiv
     a_{j}\:((q_{j},\tilde{q})) }}
 e\Bigl(\frac{\tilde{m}\tilde{a}}{\tilde{q}}\Bigr) 
     \sum_{\substack{\bar{m}<\bar{q}\\(\bar{m},\bar{q})=1\\\bar{m}\equiv
     a_{j}\:((q_{j},\bar{q})) }} 
 e\Bigl(\frac{\bar{m}\bar{a}}{\bar{q}}\Bigr) \\
  &= c_{j}(\tilde{a},\tilde{q}) \cdot c_{j}(\bar{a},\bar{q}).
\end{align*}

\vspace{-6ex}
$\hfill\square$

Proposition \ref{Prop1} provides the multiplicativity of $b$:

\begin{prop}
\label{Prop2}
Let $(\bar{q},\tilde{q})=1$. Then
$b(\bar{q}\tilde{q})=b(\bar{q})b(\tilde{q})$. 
\end{prop}
\textbf{Proof.}
We have
\begin{align*}
  b(\bar{q}\tilde{q})&= \sum_{\substack{a<\bar{q}\tilde{q}\\(a,\bar{q}\tilde{q})=1 }}
  (c_{1}c_{2}c_{3}) (a,\bar{q}\tilde{q})
  e\Bigl(-n\frac{a}{\bar{q}\tilde{q}} \Bigr) \\
  &= \sum_{\substack{\tilde{a}<\tilde{q} \\ (\tilde{a},\tilde{q})=1 }} 
     \sum_{\substack{\bar{a}<\bar{q} \\ (\bar{a},\bar{q})=1 }} 
   (c_{1}c_{2}c_{3})(\tilde{a},\tilde{q}) \cdot
   (c_{1}c_{2}c_{3})(\bar{a},\bar{q}) \cdot
    e\Bigl( -n \frac{\tilde{a}\bar{q}+\bar{a}\tilde{q}}{\bar{q}\tilde{q}} \Bigr)
\end{align*}
by substituting $a=\tilde{a}\bar{q}+\bar{a}\tilde{q}$ in the last
step. We further get
\begin{align*}
   b(\bar{q}\tilde{q})  &= 
   \sum_{\substack{\bar{a}<\bar{q} \\ (\bar{a},\bar{q})=1 }} 
    (c_{1}c_{2}c_{3})(\bar{a},\bar{q})
     e\Bigl(-n\frac{\bar{a}}{\bar{q}}\Bigr)
    \sum_{\substack{\tilde{a}<\tilde{q} \\ (\tilde{a},\tilde{q})=1 }}  
    (c_{1}c_{2}c_{3})(\tilde{a},\tilde{q})
     e\Bigl(-n\frac{\tilde{a}}{\tilde{q}}\Bigr)  \\
    &=b(\bar{q})\cdot b(\tilde{q}). 
\end{align*}

\vspace{-6ex}
\hfill$\square$

Proposition \ref{Prop2} shows that it suffices to evaluate $b$ at prime powers
$p^{k}$, $p$ prime and $k\geq 1$, to obtain formulas for $b$
and $\lambda$. It may happen that $b(p^{k})=0$, what we study now.

We first show:

\begin{prop}
\label{Prop3}
 Let $j\in\{1,2,3\}$. If $p^{k}\nmid q_{j}$ and ($p\mid q_{j}$
  or  $k\neq 1$), then $c_{j}(a,p^{k})=0$. 
\end{prop}

\textbf{Proof.}

Firstly, if $p^{k}\nmid q_{j}$ and $p\mid q_{j}$, we have 
  $d_{j}=(q_{j},p^{k})=p^{r}$ with $1\leq r<k$
  and $(d_{j},\frac{p^{k}}{d_{j}}) = (p^{r},p^{k-r})\geq p$, so $c_{j}(a,p^{k})=0.$

Secondly, if $p^{k}\nmid q_{j}$ and $k\neq 1$, then
  $d_{j}=(q_{j},p^{k})=p^{r}$ with $0\leq r<k$, and 
  \begin{equation*}
       \Big(d_{j},\frac{p^{k}}{d_{j}}\Big) = (p^{r}, p^{k-r}) = p^{\min(r,k-r)}.     
  \end{equation*}
  For $r>0$ this is $\geq p$, and so $c_{j}(a,p^{k})=0$. For
  $r=0$ we have $d_{j}=1$ and $\mu(\frac{p^{k}}{d_{j}})=\mu(p^{k})=0$
  since $k\neq 1$, so $c_{j}(a,p^{k})=0$, too.   \hfill$\square$

\zeile
Therefore $c_{j}(a,p^{k})=0$ holds unless $p^{k}\mid q_{j}$ or
($p\nmid q_{j}$ and $k=1$). This shows that
\begin{equation*}
  b(p^{k}) = \sum_{\substack{a<p^{k} \\ (a,p)=1}}
  c_{1}(a,p^{k})c_{2}(a,p^{k})c_{3}(a,p^{k}) \;e\left(-n\frac{a}{p^{k}}\right) = 0,
\end{equation*}
unless $p^{k}\mid q_{j}$ or ($p\nmid q_{j}$ and $k=1$) for every
$j=1,2,3$. We now have to consider only these cases.

\textbf{Case 1.} If $k\geq 1$, $p^{k}\mid (q_{1},q_{2},q_{3})$, then
  \begin{align*}
    b(p^{k}) &= \sum_{\substack{a<p^{k}\\(a,p)=1}}
    c_{1}(a,p^{k})c_{2}(a,p^{k})c_{3}(a,p^{k})\; 
    e\left(-n\frac{a}{p^{k}}\right) \\  &=
    \sum_{\substack{a<p^{k}\\(a,p)=1}} e\left(-n\frac{a}{p^{k}}\right) 
    \prod_{i=1,2,3} e\left(\frac{aa_{i}}{p^{k}}\right) 
     \quad(\text{since } d_{i}=(q_{i},p^{k})=p^{k}\text{ so }u_{i}=1)\\ &=
     \sum_{\substack{a<p^{k}\\(a,p)=1}}
    e\left(\frac{a_{1}+a_{2}+a_{3}-n}{p^{k}}a\right),
  \end{align*}  
  so $b(p^{k})=\ph(p^{k})$ since $p^{k}\mid a_{1}+a_{2}+a_{3}-n$
  by the general condition.

\textbf{Case 2.}
If $k=1$ and $(p,q_{1})=(p,q_{2})=(p,q_{3})=1$ then
  \begin{align*}
    b(p) &= \sum_{a=1}^{p-1} c_{1}(a,p)c_{2}(a,p)c_{3}(a,p)
    \;e\left(-n\frac{a}{p}\right) \\  &= \sum_{a=1}^{p-1} e\left(-n\frac{a}{p}\right)
    \prod_{i=1,2,3} \underbrace{\sum_{m=1}^{p-1}e\left(m\frac{a}{p}\right)}_{=-1} =
    \begin{cases}
       1-p, &\text{ if } p\mid n, \qquad \hfill(A)\\
       1, &\text{ if } p\nmid n. \qquad \hfill(B)
    \end{cases}
  \end{align*}

\textbf{Case 3.} If $k=1$, $p\mid q_{1}$ (so $d_{1}=p$) and
 $(p,q_{2})=(p,q_{3})=1$ (analogously the cases with permuted
 indices), then 
  \begin{align*}
    b(p) &=  \sum_{a=1}^{p-1} c_{1}(a,p)c_{2}(a,p)c_{3}(a,p)
    \;e\left(-n\frac{a}{p}\right) \\  &= \sum_{a=1}^{p-1}
    e\left(-n\frac{a}{p}\right)  e\left(\frac{aa_{1}}{p}\right)
     \biggl( \;\underbrace{\sum_{m=1}^{p-1}
       e\left(m\frac{a}{p}\right)}_{=-1} \;\biggr)^{2} \\ 
    &= \sum_{a=1}^{p-1} e\left(\frac{a_{1}-n}{p}a\right) =
    \begin{cases}
      p-1, &\text{ if } p\mid a_{1}-n, \qquad\hfill(C)\\
       -1, &\text{ if } p\nmid a_{1}-n. \qquad\hfill(D)
    \end{cases}
  \end{align*} 

\textbf{Case 4.} If $k=1$, $p\mid q_{1}$, $p\mid q_{2}$ and $(p,q_{3})=1$ (analogously the
  cases with permuted indices), then
  \begin{align*}
    b(p) &=  \sum_{a=1}^{p-1} c_{1}(a,p)c_{2}(a,p)c_{3}(a,p)
    \;e\left(-n\frac{a}{p}\right) \\ 
     &=  \sum_{a=1}^{p-1}
    e\left(-n\frac{a}{p}\right)e\left(\frac{aa_{1}}{p}\right)e\left(\frac{aa_{2}}{p}\right)
  \underbrace{\sum_{m=1}^{p-1} e\left(m\frac{a}{p}\right)}_{=-1}  \\
     &= -\sum_{a=1}^{p-1} e\left(\frac{a_{1}+a_{2}-n}{p}a\right) = 
     \begin{cases}
        1-p, &\text{ if } p \mid a_{1}+a_{2}-n,  \qquad\hfill (E)\\
          1, &\text{ if } p \nmid a_{1}+a_{2}-n.  \qquad\hfill (F)
     \end{cases}
    \end{align*}

If we combine all these cases, we have shown

1. If $k\geq 1$ and $p^{k}\mid (q_{1},q_{2},q_{3})$:
   $b(p^{k})=\ph(p^{k})$,
   furthermore $\lambda(p^{k})=b(p^{k})$.

2. If $p\nmid (q_{1},q_{2},q_{3})$:
\begin{equation*}
  b(p) =
  \begin{cases}
     1-p,        &\text{ if }  (p,q_{1})=(p,q_{2})=(p,q_{3})=1,\ p\mid n, \hfill(A)\\
      1,         &\text{ if }  (p,q_{1})=(p,q_{2})=(p,q_{3})=1,\ 
                             p\nmid n, \hfill(B)\\
     p-1,        &\text{ if }   p\mid q_{1},\  (p,q_{2})=(p,q_{3})=1,\ 
                              p\mid a_{1}-n, \hfill(C)\\
                 &\text{\qquad also with permuted indices,} \\
      -1,        &\text{ if }   p\mid q_{1},\  (p,q_{2})=(p,q_{3})=1,\ 
                              p\nmid a_{1}-n, \hfill(D)\\
                 &\text{\qquad also with permuted indices,} \\
      1-p,       &\text{ if }   p\mid q_{1},\  p\mid q_{2},\  (p,q_{3})=1,\ 
                              p\mid a_{1}+a_{2}-n,\; \hfill(E)\\
                 &\text{\qquad also with permuted indices,} \\
      1,         &\text{ if }   p\mid q_{1},\  p\mid q_{2},\  (p,q_{3})=1,\ 
                              p\nmid a_{1}+a_{2}-n, \;\hfill(F)\\
                 &\text{\qquad also with permuted indices,} \\
  \end{cases}
\end{equation*}
so $b(p)\in\{\plm 1,\plm(p-1)\}$. Expressed in $\lambda$ we have
\begin{equation*}
  \lambda(p)=
  \begin{cases}
    \frac{1}{(p-1)^{3}}, &(B) \\
    -\frac{1}{(p-1)^{2}}, &(A),(D) \\
    \frac{1}{p-1}, &(C),(F) \\
     -1.  &(E) \\
  \end{cases}
\end{equation*}

3. In any other case: $b(p^{k})=\lambda(p^{k})=0$.

In the following let $d:=(q_{1},q_{2},q_{3})$ where the $q_{j}$ are fixed.
For a prime $p$ let
$\gamma_{p}$ such that $p^{\gamma_{p}}\| d$, that is
$p^{\gamma_{p}}\,|\,d$ but $p^{\gamma_{p}+1}\nmid d$.

Now with $b$, $\lambda$ is multiplicative too, since
$\ph(q_{i})\ph([q_{i},\bar{q}\tilde{q}])=\ph([q_{i},\bar{q}])\ph([q_{i},\tilde{q}])$
for $(\bar{q},\tilde{q})=1$, $i=1,2,3$.
This multiplicativity for $\lambda$ shows for $Q\geq 1$:
\begin{align*}
     &\sum_{q=1}^{Q} q|\lambda(q)| 
     \leq \prod_{\substack{p\leq Q\\\text{prime}}} \left(1+
     \sum_{k=1}^{2\log Q} p^{k} |\lambda(p^{k})| \right) \\
     &= \Big(\prod_{p\leq Q, p\mid d} (1+p|\lambda(p)| + p^{2}|\lambda(p^{2})|
     +\dots + p^{\gamma_{p}}|\lambda(p^{\gamma_{p}})| )\Big)
      \cdot \Big(\prod_{p\leq Q,(p,d)=1} (1+p|\lambda(p)|) \Big)\\ 
     &\leq \Big(\prod_{p\mid d} (1+p(p-1) + p^{2}(p^{2}-p) + \dots +
     p^{\gamma_{p}}(p^{\gamma_{p}}-p^{\gamma_{p}-1})) \Big)\cdot \Big(\prod_{p\leq
     Q,(p,d)=1} (1+p|\lambda(p)|) \Big)\\ 
     &\leq \Big(\prod_{p\mid d} p^{2\gamma_{p}}\Big) \cdot \Big( \prod_{p\leq
     Q,(p,d)=1} (1+p|\lambda(p)|)  \Big)
     =d^{2}\cdot \mathcal{A}\cdot \mathcal{B}\cdot \mathcal{C}
       \cdot \mathcal{D},     
\end{align*}
where
\begin{align*}
\mathcal{A}&:= \prod_{p\leq Q,(B)} \left(1+\frac{p}{\ph(p)^{3}}\right) \leq
  \sum_{q,p\mid q\Rightarrow p\leq Q} \frac{q\mu^{2}(q)}{\ph(q)^{3}} \ll
  \sum_{q} \frac{\mu^{2}(q)}{q^{2}} (\loglog q)^{3} \ll 1, \\
\mathcal{B}&:= \prod_{p\leq Q,(A)} \left(1+\frac{p}{\ph(p)^{2}}\right)
   \cdot\prod_{i=1,2,3} \prod_{\substack{p\leq Q,\\ (D)\text{ for }
  q_{i}}} \left(1+\frac{p}{\ph(p)^{2}}\right)\\
     &\leq \left( \sum_{\substack{q\leq n\\p\mid q\Rightarrow p\mid n}}
  \frac{q\mu^{2}(q)}{\ph(q)^{2}} \right) \prod_{i=1,2,3}
  \left( \sum_{\substack{q\leq q_{i}\\p\mid q\Rightarrow p\mid q_{i}}}
  \frac{q\mu^{2}(q)}{\ph(q)^{2}} \right) 
     \ll \left(\sum_{q\leq n} \frac{\mu^{2}(q)}{q}(\loglog q)^{2}
  \right)^{4}  \\ &\ll (\log n)^{8}, \\
\mathcal{C}&:=\prod_{i=1,2,3} \prod_{\substack{p\leq Q\\ (C)\text{ for }
  q_{i}}} \left(1+\frac{p}{p-1}\right) \cdot
  \prod_{\substack{i,j\in\{1,2,3\}\\ i\neq j}}
  \prod_{\substack{p\leq Q\\(F)\text{ for } q_{i},q_{j}}}
  \left(1+\frac{p}{p-1} \right) \\
   &\leq \prod_{i=1,2,3} \prod_{p\mid q_{i}} (1+2) \cdot
  \prod_{\substack{i,j\in\{1,2,3\}\\ i\neq j}} 
   \prod_{p\mid (q_{i},q_{j})} (1+2) 
   \leq\prod_{i=1,2,3} 2^{2\omega(q_{i})} \cdot
  \prod_{\substack{i,j\in\{1,2,3\}\\ i\neq j}} 
  2^{2\omega((q_{i},q_{j}))} \\
   &\leq \tau^{2}(q_{1}) \tau^{2}(q_{2}) \tau^{2}(q_{3}) \tau^{2}((q_{1},q_{2}))
  \tau^{2}((q_{1},q_{3})) \tau^{2}((q_{2},q_{3})) \leq
  \tau^{4}(q_{1})\tau^{4}(q_{2})\tau^{4}(q_{3}),
 \\
\mathcal{D}&:= \prod_{\substack{i,j,k\in\{1,2,3\}\\ i,j,k \text{
      p.w.d.}}} 
  \prod_{\substack{p\leq Q\\(E)\text{ for } q_{i},q_{j},q_{k}}} (1+p)
  \leq  \prod_{\substack{i,j,k\in\{1,2,3\}\\ i,j,k \text{ p.w.d.}}} 
  \prod_{\substack{p\mid (q_{i},q_{j})\\ p\nmid q_{k}}} (1+p) \\
  &= \prod_{\substack{i,j,k\in\{1,2,3\}\\ i,j,k \text{ p.w.d.}}} 
   \sigma\biggl( \prod_{\substack{p\mid (q_{i},q_{j}) \\ p\nmid q_{k}}}  p
  \biggr) \leq \prod_{\substack{i,j,k\in\{1,2,3\}\\ i,j,k \text{ p.w.d.}}} 
  \sigma\left(\frac{(q_{i},q_{j})}{(q_{1},q_{2},q_{3})}
  \right) \ll \prod_{\substack{i,j,k\in\{1,2,3\}\\ i,j,k
  \text{ p.w.d.}}} \frac{(q_{i},q_{j})}{d}\log n \\
  &= \frac{1}{d^{3}}(q_{1},q_{2})(q_{1},q_{3})(q_{2},q_{3}) (\log n)^{3},
\end{align*}
where $\sigma(t):=\sum_{t\mid q} t$ is the divisor sum function, for which
$\sigma(t)\ll t\log t$ holds, and $\omega(t)$ is the number of distinct prime
factors of $t$.

Therefore
\begin{equation*}
  \sum_{q=1}^{Q} q|\lambda(q)| \ll  
  \frac{(q_{1},q_{2})(q_{1},q_{3})(q_{2},q_{3})}{(q_{1},q_{2},q_{3})}
  \tau^{4}(q_{1})\tau^{4}(q_{2})\tau^{4}(q_{3}) (\log n)^{11},
\end{equation*}
also true for $Q\to\infty$. So for any $Q\geq 1$ we have
\begin{equation*}
  \sum_{q\geq Q} |\lambda(q)|\leq \frac{1}{Q} \sum_{q=1}^{\infty}
  q|\lambda(q)| \leq
  \frac{1}{Q}\frac{(q_{1},q_{2})(q_{1},q_{3})(q_{2},q_{3})}{(q_{1},q_{2},q_{3})}  
  \tau^{4}(q_{1})\tau^{4}(q_{2})\tau^{4}(q_{3}) (\log n)^{11}.
\end{equation*}

We see that the singular series
$\mathcal{S}_{3}(n)=\sum_{q=1}^{\infty}\lambda(q)$ converges
absolutely, and we have
\[
    \mathcal{S}_{3}(n) \ll
    \frac{(q_{1},q_{2})(q_{1},q_{3})(q_{2},q_{3})}{(q_{1},q_{2},q_{3})}   
  \tau^{4}(q_{1})\tau^{4}(q_{2})\tau^{4}(q_{3}) (\log n)^{11}.
\]

It follows further that
\begin{align*}%
  \sum_{q_{1},q_{2},q_{3}} \max_{a_{1},a_{2},a_{3}} e_{1}
   &\ll \sum_{q_{1},q_{2},q_{3}}
  \frac{(q_{1},q_{2})(q_{1},q_{3})(q_{2},q_{3})}{d\ph(q_{1})\ph(q_{2})\ph(q_{3})}
  \frac{n^{2}}{R}\tau^{4}(q_{1})\tau^{4}(q_{2})\tau^{4}(q_{3}) (\log n)^{11} \\
  &\ll \frac{n^{2}}{R} (\log n)^{12} \sum_{q_{1},q_{2},q_{3}}
   \frac{\tau^{4}(q_{1})\tau^{4}(q_{2})\tau^{4}(q_{3})
   (q_{1},q_{2})(q_{1},q_{3})(q_{2},q_{3})}{q_{1}q_{2}q_{3}d} 
   \\ &\ll \frac{n^{2}}{R}(\log n)^{\eta} \ll \frac{n^{2}}{(\log n)^{A+3}},  
\end{align*}
since $B\geq A+\eta+3$ in $R=(\log n)^{B}$ for some absolute constant
$\eta>0$. This can be proven as follows. By using
\begin{equation*}
  \sum_{t\leq n} \frac{\tau^{m}(t)}{t} \leq (\log n)^{2^{m}}, 
\end{equation*}
we see that
\begin{align*}
  &\sum_{q_{1},q_{2},q_{3}}
  \frac{(q_{1},q_{2})(q_{1},q_{3})(q_{2},q_{3})}
  {q_{1}q_{2}q_{3}(q_{1},q_{2},q_{3})}
  \tau^{4}(q_{1})\tau^{4}(q_{2})\tau^{4}(q_{3}) \\
  &\leq \sum_{d\leq n} \sum_{a,b,c\leq n} \sum_{e,f,g\leq n}
  \frac{dadbdc}{d^{3}a^{2}b^{2}c^{2}efgd}
  \tau^{12}(d)\tau^{8}(a)\tau^{8}(b)\tau^{8}(c)\tau^{4}(e)\tau^{4}(f)\tau^{4}(g) \\
  &= \sum_{d}\sum_{a,b,c} \sum_{e,f,g} \frac{
  \tau^{12}(d)\tau^{8}(a)\tau^{8}(b)\tau^{8}(c)\tau^{4}(e)\tau^{4}(f)\tau^{4}(g) }
  {abcdefg} \\ &\ll (\log n)^{\eta}
\end{align*}
for some absolute constant $\eta>0$,
where we substituted $q_{1}=dabe$, $q_{2}=dacf$, $q_{3}=dcbg$ with pairwise
coprime $a,b,c$ and $e,f,g$.

Further we have
\begin{equation*}
  \sum_{q\leq R} q^{2}|\lambda(q)| \leq R \sum_{q=1}^{\infty}
  q|\lambda(q)| \ll R \tau^{4}(q_{1})\tau^{4}(q_{2})\tau^{4}(q_{3})
  \frac{(q_{1},q_{2})(q_{1},q_{3})(q_{2},q_{3})}{d}(\log n)^{11},
\end{equation*}
so also
\begin{align*}
  \sum_{q_{1},q_{2},q_{3}} \max_{a_{1},a_{2},a_{3}} e_{2}
   &\ll \sum_{q_{1},q_{2},q_{3}}
  \frac{n^{2}(q_{1},q_{2})(q_{1},q_{3})(q_{2},q_{3})}
  {R \ph(q_{1})\ph(q_{2})\ph(q_{3})d}
  \tau^{4}(q_{1})\tau^{4}(q_{2})\tau^{4}(q_{3}) (\log n)^{11} \\ &\ll
  \frac{n^{2}}{(\log n)^{A+3}}
\end{align*}
as above.

So everything concerning Theorem \ref{Th2} is shown. \hfill$\square$
 
\subsubsection{Discussion of the singular series}
\label{singseries}

Now we consider $\mathcal{S}_{3}(n)$ under the general condition.

Since $\mathcal{S}_{3}(n)$ is absolutely convergent and since $\lambda$ is
multiplicative, we see that it has an Eulerproduct, namely
\begin{equation*}
  \mathcal{S}_{3}(n)=\prod_{p}\Bigl( 1+\sum_{k=1}^{\infty} \lambda(p^{k}) \Bigr).
\end{equation*}
For $p^{\alpha}\| (q_{1},q_{2},q_{3})$ we have
$1+\lambda(p)+\dots+\lambda(p^{\alpha}) = p^{\alpha}$ and for other
primes $p$ we get factors according to the
cases $(A),\dots,(F)$. Moreover we see that
$\mathcal{S}_{3}(n)$ vanishes if case $(E)$ for a prime $p$ occurs, that
is if

\begin{align*}
  (E): \qquad\exists\; j,k,l\in\{1,2,3\} \text{ pairwise different with } \\
      p\mid (q_{j},q_{k}),\: p\nmid q_{l},\: p\mid n-(a_{j}+a_{k}).
\end{align*}

In all other cases we have
\begin{equation*}
  \mathcal{S}_{3}(n)= 
  (q_{1},q_{2},q_{3})
 \prod_{\substack{p, (A) \\ \text{ or } (D)}} \Bigl(1-\frac{1}{(p-1)^{2}}\Bigr)
 \prod_{p, (B)} \Bigl(1+\frac{1}{(p-1)^{3}}\Bigr)
  \prod_{\substack{p, (C) \\ \text{ or } (F)}} \Bigl(1+\frac{1}{p-1}\Bigr)
\end{equation*}
with properties
\begin{align*}
  &(A): (p,q_{1})=(p,q_{2})=(p,q_{3})=1, p\mid n,\\
  &(B): (p,q_{1})=(p,q_{2})=(p,q_{3})=1, p\nmid n,\\
  &(C): \exists\; j,k,l\in\{1,2,3\} \text{ pwd: } 
      p\mid q_{j},\: (p,q_{k})=(p,q_{l})=1, \: p\mid n-a_{j}, \\
  &(D): \exists\; j,k,l\in\{1,2,3\} \text{ pwd: } 
      p\mid q_{j},\: (p,q_{k})=(p,q_{l})=1, \: p\nmid n-a_{j}, \\
  &(F): \exists\; j,k,l\in\{1,2,3\} \text{ pwd: } 
      p\mid q_{j},\: p\mid q_{k},\: (p,q_{l})=1,\: p\nmid n-(a_{j}+a_{k}). 
\end{align*}
So we see that $\mathcal{S}_{3}(n)=0$ if and only if case $(E)$ occurs
or the general condition is not fulfilled.
Further if $\mathcal{S}_{3}(n)>0$ we see from the Eulerproduct that it
is at least some absolute positive constant times
$(q_{1},q_{2},q_{3})$, since $\prod_{p>2} (1-(p-1)^{-2})$ converges
and the other products are $>1$.

Now we prove

\textbf{Lemma \ref{lemstar}.}\!\!
\textit{
  If $n$ is odd, then for given $q_{3}$, $a_{3}$ with
  $(a_{3},q_{3})=1$ and $q_{2}$ there exists an admissible $a_{2}$
  (such that for every $q_{1}$ there exists an admissible $a_{1}$).
  For even $n$ and given $q_{1},q_{2},q_{3}$ 
  there exists \emph{no} admissible triplet $a_{1},a_{2},a_{3}$.
}

Recall that $a_{1},a_{2},a_{3}$ is admissible for $q_{1},q_{2},q_{3}$,
if $(a_{i},q_{i})=1$ for $i=1,2,3$, $n\equiv a_{1}+a_{2}+a_{3} \mo
(q_{1},q_{2},q_{3})$ and $\mathcal{S}_{3}(n)>0$.

\textbf{Proof.}
For the proof, let $q:=(q_{1},q_{2},q_{3})$ and denote by
$\nu_{p}(m)$ the exponent of a prime $p$ in $m$, that is
$p^{\nu_{p}(m)}\mid m$ but $p^{\nu_{p}(m)+1}\nmid m$.

First let $n$ be even, and consider $q_{1},q_{2},q_{3}$ with\\[-5ex]
\begin{enumerate}[(a)]
\item $2\mid q_{j}, 2\nmid q_{k},q_{l}$. Then $(A),(B),(E),(F)$ are
  not possible, and the condition $2\mid n-a_{j}$ is wrong since
  $a_{j}$ must be odd. Therefore $(D)$ holds with $p=2$, and so
  $\mathcal{S}_{3}(n)=0$.\\[-4ex]
\item
 $2\mid q_{j},q_{k}$ and $2\nmid q_{l}$. Then $(A),\dots,(D)$ are
  not possible, and condition $2\mid n-(a_{j}+a_{k})$ in $(E)$
  holds since $a_{j},a_{k}$ are odd, so $\mathcal{S}_{3}(n)=0$.\\[-4ex]
\item
 Further $2\mid q_{1},q_{2},q_{3}$ is not possible since then
  $a_{1},a_{2},a_{3}$ are odd and so $n\not\equiv
  a_{1}+a_{2}+a_{3}\:(q)$, so  
  $\mathcal{S}_{3}(n)=0$.\\[-4ex]
\item
 Also $2\nmid q_{1},q_{2},q_{3}$ is not possible since then $(A)$
  holds for $p=2$, so $\mathcal{S}_{3}(n)=0$ holds.
\end{enumerate}

Now let $n$ be odd and let $q_{3},a_{3}$ with $(a_{3},q_{3})=1$ and
$q_{2}$ be given. We construct $a_{2}$ and $q_{2}$ with
$(a_{2},q_{2})=1$ such that
\begin{equation*}
  \forall\, p\mid (q_{3},q_{2}): n\not\equiv a_{3}+a_{2}\ (p).
\end{equation*}
For any $p\mid(q_{3},q_{2})$ take $h_{p}$ such that $1\leq h_{p}\leq p-1$
with $n-a_{3}+h_{p}\not\equiv 0\ (p)$.
Such a number $h_{p}$ exists for $p>2$ since then $p-1>1$,
and if $p=2$ take $h_{2}=1$ since $n-a_{3}+1\not\equiv 0\:(2)$ holds for $p=2\mid
(q_{3},q_{2})$, where $q_{3}$ is even and therefore $a_{3}$ is odd.

Then take $a_{2}$ with $(a_{2},q_{2})=1$ and $a_{2}\equiv
n-a_{3}+h_{p}\ (p)$ for every $p\mid(q_{3},q_{2})$ via the Chinese
Remainder Theorem. Now we prove that this $a_{2}$ is admissible. For
this, consider now any $q_{1}$, and we have to find now an admissible
$a_{1}$, that means such that

\qquad $(1)$ $n\equiv a_{1}+a_{2}+a_{3}\: ((q_{1},q_{2},q_{3}))$,

\qquad $(2)$ $\forall\, p\mid (q_{1},q_{2}), \,p\nmid q_{3}:
n\not\equiv a_{1}+a_{2} \:(p)$, 

\qquad $(3)$ $\forall\, p\mid (q_{1},q_{3}), \,p\nmid q_{2}:
n\not\equiv a_{1}+a_{3}  \:(p)$, 

\qquad $(4)$ $\forall\, p\mid (q_{2},q_{3}), \,p\nmid q_{1}:
n\not\equiv a_{2}+a_{3}  \:(p)$.

Now condition $(4)$ is fulfilled by the choice of $a_{2}$. We have
to construct now an admissible $a_{1}$ mod $q_{1}$, $(a_{1},q_{1})=1$,
namely such that conditions $(1)-(3)$ are fulfilled.

Firstly, $a_{1}$ has to be such that $a_{1}\equiv n-a_{2}-a_{3} 
\ ((q_{1},q_{2},q_{3}))$. Since $n-a_{2}-a_{3}\equiv -h_{p}\not\equiv
0\ (p)$ for any $p\mid(q_{2},q_{3})$ we see that $a_{1}$ mod
$(q_{1},q_{2},q_{3})$ may be chosen like that, and it will not
contradict to $(a_{1},q_{1})=1$, and also condition $(1)$ is fulfilled.

Further $a_{1}$ must be $a_{1}\equiv n-a_{3}+k_{p}\not\equiv 0\ (p)$ 
for every $p\mid (q_{1},q_{3})$, $p\nmid q_{2}$, where $1\leq
k_{p}\leq p-1$ (condition $(3)$),
and also with $a_{1}\equiv n-a_{2}+l_{p}\not\equiv 0\ (p)$ for every
$p\mid (q_{2},q_{1})$, $p\nmid q_{3}$, where $1\leq l_{p}\leq p-1$
(condition $(2)$). Here the existence of $l_{p}$ and $k_{p}$
can be explained as above for $h_{p}$. Then take $a_{1}$ with
$(a_{1},q_{1})=1$ to hold these congruences, again via the Chinese
Remainder Theorem. It is admissible by construction.
\hfill$\square$

\zeile
By studying property $(E)$, we encounter the following connection with
the binary Goldbach problem. 

Let $p$ be any prime $>2$ and let $n$ be sufficiently large.
We can construct $a_{i},q_{i}$, with
$(a_{i},q_{i})=1$ for $i=1,2,3$, and with
\begin{equation*}
  p\mid(q_{1},q_{2}), \enspace p\nmid q_{3}, \enspace n\equiv a_{1}+a_{2}\;(p),
  \enspace a_{1}+a_{2}+a_{3} \equiv  n \;((q_{1},q_{2},q_{3})),
\end{equation*}
namely take any odd $q_{1},q_{2},q_{3}$ such that $p\mid(q_{1},q_{2}), \enspace p\nmid q_{3}$,
and take $a_{1}$ with $n-a_{1}\not\equiv 0\:(p)$ relatively prime to
$q_{1}$, take $a_{2}$ relatively prime to $q_{2}$ with $a_{2}\equiv
n-a_{1} \:(p)$ and $(n-a_{1}-a_{2},(q_{1},q_{2},q_{3}))=1$, and $a_{3}$ with 
$a_{3}\equiv n-a_{1}-a_{2} \;((q_{1},q_{2},q_{3}))$ relatively prime to $q_{3}$.
If we could show that there exist primes
$p_{i}\equiv a_{i}\;(q_{i})$, $i=1,2,3$, with
$n=p_{1}+p_{2}+p_{3}$, and so $n\equiv
a_{1}+a_{2}+p_{3}\;((q_{1},q_{2}))$, then since $n\equiv a_{1}+a_{2}\;(p)$
it follows that $0\equiv p_{3}\;(p)$, so $p_{3}=p$ and
$n-p=p_{1}+p_{2}$. Then the number $n-p$ would be the sum 
of two primes.

So if the considered ternary Goldbach problem with primes in
independent arithmetic progressions touches the binary Goldbach
problem, the circle method fails.

%
\subsection{A Lemma involving sieve methods}
\label{secthree}

Before considering the minor arcs we show the following Lemma
by using the large sieve inequality and a
formula of Montgomery in \cite{c7}.
The method was already presented in \cite{c2}.

\begin{lemma}
\label{lem1}
  For $Q\geq 1$, $H>0$ and $b_{1},\dots,b_{n}\in\C$ we have
  \label{eq:eqth1}
  \begin{align*}
    \sum_{q\sim Q} &q \max_{0\leq a<q}\;
    \Biggl|\sum_{\substack{m\leq n \\ m\equiv a (q)}}
    b_{m}\Biggr|^{\,2} \\
    &\ll (n^{2}+Q^{2})H^{-1}(\log Q) \max_{m\leq n} |b_{m}|^{2}
     + (n+Q^{2})\, H (\log Q) \sum_{m\leq n} |b_{m}|^{2}  
  \end{align*}
  \end{lemma}

with an absolute $O$-constant.

\zeile
\textbf{Remark.} If $Q$ may be some small power of $n$ the Cauchy-Schwarz-estimate
\begin{equation*} 
   \sum_{q\sim Q} q \max_{0\leq a<q} \;
    \Biggl|\sum_{\substack{m\leq n \\ m\equiv a (q)}}
    b_{m}\Biggr|^{\,2} \ll \sum_{q\leq 2Q} q \sum_{m\leq n} 
    |b_{m}|^{2} \; \frac{n}{q} \ll  n\, Q \sum_{m\leq n} |b_{m}|^{2}
\end{equation*}
is weaker. An approach with the large sieve inequality involving
characters does not work either.

\zeile\zeile 
\textbf{Proof of Lemma \ref{lem1}.} 

For a residue class $a$ mod $q$ we set
\[
   N(a,q) := \sum_{\substack{m\leq n\\m\equiv a (q)}} b_{m}.
\]

Now the expression on the left hand side in
Lemma \ref{lem1} is $E_{1}+E_{2}$ with
\begin{equation*}
  E_{1}:= \sum_{\substack{q\sim Q\\d(q)>H}}
    q \max_{0\leq a<q}  |N(a,q)|^{2}
\end{equation*}
and
\begin{equation*}
  E_{2}:= \sum_{\substack{q\sim Q\\d(q)\leq H}}
    q \max_{0\leq a<q}  |N(a,q)|^{2}.
\end{equation*}

Consider first $E_{1}$. Let
\begin{equation*}
  A_{Q}:=\#\{ q \; ; \; q\sim Q, \;d(q) > H \},
\end{equation*}
then
\begin{equation*}
   A_{Q} H <
  \sum_{\substack{q\sim Q \\ d(q)> H}} d(q)
  \leq\sum_{q\leq 2Q} d(q) \ll  Q\log Q,
\end{equation*}
so
\begin{equation*}
  A_{Q}\ll \frac{Q\log Q}{H}.
\end{equation*}
Since $N(a,q)\ll(\frac{n}{q}+1)\max_{m\leq n} |b_{m}|$ we get
\begin{align*}
  E_{1}&\ll \sum_{\substack{q \sim Q \\ d(q)>H }}
      q  \max_{a} |N(a,q)|^{2} \ll
      \sum_{\substack{q \sim Q \\ d(q)>H }} q\biggl(\frac{n^{2}}{q^{2}}
      +1\biggr) \max_{m\leq n} |b_{m}|^{2}  \\ &\ll
      A_{Q} \biggl(\frac{n^{2}}{Q} +Q\biggr)\max_{m\leq n} |b_{m}|^{2}
      \ll \biggl(\frac{n^{2}}{H} +\frac{Q^{2}}{H}\biggr) (\log Q) \max_{m\leq n}
 |b_{m}|^{2} .
\end{align*}
This is the first summand on the right hand side of Lemma \ref{lem1}.

\zeile\zeile
Now to $E_{2}$.

For any integer $0\leq h<q$ let
\[
   f_{h}(q):=\sum_{d|q} \mu(d) \frac{q}{d} N\Bigl(h,\frac{q}{d}\Bigr),
\]
so M\"obius' inversion formula gives
\begin{equation*}
  q N(h,q) = \sum_{d|q} f_{h}(d)
\end{equation*}
for all $0\leq h<q$.
With this we have
\begin{align*}
  \label{eq:eq1}
  E_{2}&=\sum_{\substack{q\sim Q\\ d(q)\leq H}}
    \frac{1}{q} \max_{0\leq a<q} q^{2} |N(a,q)|^{2}
  =\sum_{\substack{q\sim Q\\ d(q)\leq H}}
   \frac{1}{q} \max_{0\leq a<q} \;\Biggl| \sum_{d|q} f_{a}(d)\Biggr|^{\,2} \\
  &\leq \sum_{\substack{q\sim Q\\
      d(q)\leq H}}  \frac{d(q)}{q}  \sum_{d|q} \max_{0\leq a<q} |f_{a}(d)|^{2}.
\end{align*}

The maximum is taken over $a$ with $0\leq a<q$.
We see that $|f_{a}(d)|^{2}$ is $d$-periodic in $a$ for $d|q$,
since  $N(a+d,t)=N(a,t)$ for $t|d$, so
\begin{equation*}
  f_{a+dl}(d)=\sum_{t|d}\mu(t) \frac{d}{t} N\Bigl(a+dl,\frac{d}{t}\Bigr)
  = \sum_{t|d}\mu(t) \frac{d}{t} N\Bigl(a,\frac{d}{t}\Bigr)
  = f_{a}(d)\text{ for all } l\in\Z,
\end{equation*}
therefore the maximum stays equal if taken only over $a$ with
$0\leq a<d$. We estimate this maximum by $\sum_{0\leq a<d}$ and get
\begin{align*}
  E_{2} &\leq \sum_{\substack{q\sim Q\\ d(q)\leq H}} \frac{d(q)}{q}
     \,  \sum_{d|q} \sum_{0\leq a<d} |f_{a}(d)|^{\,2}.
\end{align*}

By Montgomery in \cite{c7}, equation (10), we have 
for $T(\alpha):=\sum_{m\leq n} b_{m}e(\alpha m)$, $\alpha\in\R$, the formula
\begin{equation*}
  \frac{1}{d} \sum_{h=0}^{d-1} |f_{h}(d)|^{2}
    = \sum_{\substack{a<d \\ (a,d) = 1}}
    \biggl|T\biggl(\frac{a}{d}\biggr)\biggr|^{2},
\end{equation*}
that we can apply here. We get
\begin{align*}
   E_{2}&\leq \sum_{\substack{q\sim Q\\ d(q)\leq H}} d(q)
   \;\sum_{d|q} \frac{d}{q} \sum_{\substack{a<d \\ (a,d) = 1}}
     \biggl|T\biggl(\frac{a}{d}\biggr)\biggr|^{2} \\
   &\leq H\sum_{d\leq 2Q} \biggl(\sum_{\substack{q\sim Q\\
    d|q}} \frac{d}{q} \biggr)
        \sum_{\substack{a<d \\ (a,d) = 1}}
         \biggl|T\biggl(\frac{a}{d}\biggr)\biggr|^{2} \\
   &\ll H (\log Q)\sum_{d\leq 2Q}
        \sum_{\substack{a<d \\ (a,d) = 1}}
         \biggl|T\biggl(\frac{a}{d}\biggr)\biggr|^{2} \\
   &\ll H(\log Q) \; (n+Q^{2}) \sum_{m\leq n} |b_{m}|^{2}
\end{align*}
by the inequality of the large sieve. This is the second term on
the right hand side of Lemma \ref{lem1}. \hfill$\square$

%
\subsection{The conclusion with Lemma \ref{lem1}}
\label{secfour}

Now let $A,\theta >0$  and $\vartheta>0$ as above.
Let $Q_{1},Q_{2},Q_{3}\leq n^{\h}/(\log n)^{\vartheta}$.

We consider first
\begin{equation*}
  \mathcal{E}_{Q_{1},Q_{2},Q_{3}}^{\mathfrak{m}}:=
  \sum_{q_{3}\sim Q_{3}} \max_{a_{3}}\sum_{q_{2}\sim Q_{2}}
  \max_{a_{2}} \sum_{q_{1}\sim Q_{1}} \max_{a_{1}}
    \left| J_{3}^{\mathfrak{m}}(n) \right|.
\end{equation*}
From the definition of $J_{3}$ and $J_{2}$ we have
\begin{align*}
  \mathcal{E}_{Q_{1},Q_{2},Q_{3}}^{\mathfrak{m}} &\leq
\sum_{q_{3}} \max_{a_{3}}\sum_{q_{2}}
  \max_{a_{2}} \sum_{q_{1}} \max_{a_{1}}
   \sum_{\substack{m_{1}\leq n \\ m_{1}\equiv a_{1} (q_{1})}} \Lambda(m_{1}) 
          \left| J_{2}^{\mathfrak{m}}(n-m_{1})
  \right| \\
   &\leq \sum_{q_{3}} \max_{a_{3}}\sum_{q_{2}} \max_{a_{2}}
  \sum_{q_{1}} \max_{a_{1}} \sum_{\substack{m\leq n \\ m\equiv
    n-a_{1} (q_{1})}} (\log n)\; |J_{2}^{\mathfrak{m}}(m)|.
\end{align*}

By Cauchy-Schwarz' inequality we now get
\begin{equation*}
  \mathcal{E}_{Q_{1},Q_{2},Q_{3}}^{\mathfrak{m}} \leq (\log n)
   \sum_{q_{3}} \max_{a_{3}}\sum_{q_{2}} \max_{a_{2}}
    \biggl( \sum_{q_{1}\sim Q_{1}} q_{1} \max_{a_{1}}
  \biggl|\sum_{\substack{m\leq n \\ m\equiv a_{1}\;(q_{1})}}
   |J_{2}^{\mathfrak{m}}(m)| \biggr|^{2}  \: \biggr)^{\h} 
\end{equation*}
and we apply Lemma \ref{lem1} to the expression in large brackets.

Since $Q_{1}\leq n^{\h}$ we see that
\begin{align*}
 \mathcal{E}_{Q_{1},Q_{2},Q_{3}}^{\mathfrak{m}} \ll (\log n)
  \sum_{q_{3}} \max_{a_{3}}\sum_{q_{2}} \max_{a_{2}}
 \biggl( &\frac{n^{2}}{H}(\log n)\max_{m\leq n} |J_{2}^{\mathfrak{m}}(m)|^{2} \\
    &+ nH (\log n) \sum_{m\leq n} |J_{2}^{\mathfrak{m}}(m)|^{2} \biggr)^{\h}.
\end{align*}

Now we apply the following two lemmas, which will be proven
in the last paragraphs.

\begin{lemma}
  \label{lem2}
   For $Q_{2},Q_{3}\leq n^{\h}/(\log n)^{\vartheta}$ we have
  \begin{equation*}
    \sum_{q_{2},q_{3}}\max_{\substack{m\leq n\\a_{2},a_{3}}} 
  |J_{2}^{\mathfrak{m}}(m)|\ll n(\log n)^{7}.
  \end{equation*}
\end{lemma}

\begin{lemma}
  \label{lem3}
  For $Q_{2}\leq n^{\h}/(\log n)^{\vartheta}$ and $Q_{3}\leq (\log
  n)^{\theta}$ we have
  \begin{equation*}
    \sum_{q_{3}} \max_{a_{3}} \sum_{q_{2}} \max_{a_{2}}\biggl(
     \sum_{m\leq n} |J_{2}^{\mathfrak{m}}(m)|^{2}\biggr)^{\!\hb}
    \ll \frac{n^{3/2}}{(\log n)^{2A+16}}.
  \end{equation*}
\end{lemma}
Here the sum over such a small $Q_{3}$-range is of course pointless;
but we state it here to see why no larger bound for $Q_{3}$ is
possible to get with the given method in the proof of Lemma 4.

With $H:=(\log n)^{2A+23}$
it follows from Lemma \ref{lem2} and \ref{lem3} that
\begin{equation*}
  \mathcal{E}_{Q_{1},Q_{2},Q_{3}}^{\mathfrak{m}} \ll
  \frac{n^{2}}{(\log n)^{A+3}}.
\end{equation*}

Finally, together with Theorem \ref{Th2}, 
we get for $Q_{2}\leq n^{\h}/(\log n)^{\vartheta}$ 
and $Q_{3}\leq (\log n)^{\theta}$ the estimate
\begin{align*}
  \sum_{q_{3}\sim Q_{3}} \; \max_{\substack{a_{3}\\(a_{3},q_{3})=1}} \;
   &\sum_{q_{2}\sim Q_{2}} \; \max_{\substack{a_{2}\\(a_{2},q_{2})=1}} \;
  \sum_{q_{1}\sim Q_{1}} \; \max_{\substack{a_{1}\\(a_{1},q_{1})=1}} 
          \left|  J_{3}(n)   - \frac{n^{2}\mathcal{S}_{3}(n)}
  {2\ph(q_{1})\ph(q_{2})\ph(q_{3})} \right| \\
  &\leq \sum_{\substack{k,Q_{3}=2^{k}\\\leq (\log n)^{\theta}}}
   \:\sum_{\substack{j,Q_{2}=2^{j}\\\leq  n^{\h}/(\log n)^{\vartheta}}}
 \: \sum_{\substack{i,Q_{1}=2^{i}\\\leq n^{\h}/(\log n)^{\vartheta}}}
   \left(\mathcal{E}_{Q_{1},Q_{2},Q_{3}}^{\mathfrak{m}}
    + \mathcal{E}_{Q_{1},Q_{2},Q_{3}}^{\mathfrak{M}} \right) \\
 &\ll (\log n)^{3}\cdot\frac{n^{2}}{(\log n)^{A+3}} =
   \frac{n^{2}}{(\log n)^{A}},
\end{align*}
and from that follows Theorem \ref{Th1}.

So it remains to show Lemma \ref{lem2} and Lemma \ref{lem3}.

%
%
\subsection{Two Lemmas on the minor arcs}
\label{secfive}
\subsubsection{Proof of Lemma \ref{lem2}}

We have
\begin{equation*}
   \sum_{q_{2},q_{3}} \max_{\substack{m\leq n \\ a_{2},a_{3}}}
   |J_{2}^{\mathfrak{m}} (m)| = \sum_{q_{2},q_{3}}
   \max_{\substack{m\leq n \\ a_{2},a_{3}}}
   |J_{2}(m) -  J_{2}^{\mathfrak{M}} (m)|.
\end{equation*}

Now we estimate $J_{2}(m)$ and $J_{2}^{\mathfrak{M}}(m)$.
The reason why we split $J_{2}^{\mathfrak{m}}(m)$ is that the
trivial upper estimate for $J_{2}^{\mathfrak{m}}(m)$, namely 
\begin{equation*}
  J_{2}^{\mathfrak{m}}(m) \ll \int_{0}^{1}
     |S_{2}(\alpha)S_{3}(\alpha)| d\alpha,
\end{equation*}
does not suffice.

We have
\begin{align*}
   J_{2}(m) &= \int_{0}^{1} S_{2}(\alpha)S_{3}(\alpha)\; e(-m\alpha)
   \,d\alpha \\
  &= \sum_{\substack{ m_{2}\leq n \\ m_{2}\equiv
       a_{2} \:(q_{2})}} \Lambda(m_{2}) \sum_{\substack{ m_{3}\leq n \\
       m_{3}\equiv a_{3} \:(q_{3}) }} \Lambda(m_{3}) \int_{0}^{1}
   e(\alpha(m_{2}+m_{3}-m))\, d\alpha
\end{align*}
and by the orthogonal relations for $e(\alpha m)$ we have that the
last integral is $1$, if $m_{2}+m_{3}=m$, and $0$ otherwise. Therefore we
get
\[
  J_{2}(m) \ll \sum_{\substack{ m_{2}\leq n \\ m_{2}\equiv
      a_{2}\:(q_{2}) \\ m_{2} \equiv m-a_{3} \:(q_{3}) }} (\log n)^{2} 
   \ll \frac{n}{[q_{2},q_{3}]} (\log n)^{2} \ll
   \frac{n}{q_{2}q_{3}}(\log n)^{2} (q_{2},q_{3}),
\]
so
\begin{align*}
  &\sum_{q_{2},q_{3}} \max_{\substack{m\leq n \\a_{2},a_{3}}} |J_{2}(m)| 
   \ll n (\log n)^{2} \sum_{q_{2},q_{3}}
   \frac{(q_{2},q_{3})}{q_{2}q_{3}} \\ &\ll n (\log n)^{2}
   \sum_{m,q_{2}',q_{3}'} \frac{m}{mq_{2}'mq_{3}'} \ll n(\log n)^{5}.
\end{align*}

Now we consider the following
\begin{prop}
  We have
  \[  \sum_{q_{2},q_{3}} \max_{\substack{a_{2},a_{3}\\m\leq n}}
   |J_{2}^{\mathfrak{M}} (m)| \ll n(\log n)^{7}. 
  \]
\end{prop}
By this and together with above estimation we get therefore
Lemma \ref{lem2}.\hfill$\square$

\textbf{Proof of Proposition 4.}

We have to consider the analogous estimation for
$J_{2}^{\mathfrak{M}}(m)$ as was done in paragraph \ref{parzwo} in order
to estimate $J_{3}^{\mathfrak{M}}(m)$. 

We get
\[
  J_{2}^{\mathfrak{M}} (m) = \sum_{q\leq R} \sum_{\substack{a<q \\
    (a,q)=1}} I(a,q)
\]
with
\begin{align*}
  I(a,q) &= \int_{-R/qn}^{R/qn} S_{2} \biggl(\frac{a}{q} + \alpha\biggr) 
  S_{3} \biggl(\frac{a}{q} + \alpha\biggr) e\biggl(-m\biggl(\frac{a}{q} +
  \alpha\biggr)\biggr)d\alpha \\
   &= \frac{(c_{2}c_{3})(a,q)}{\ph([q_{2},q])\ph([q_{3},q])}
   e\biggl(-m\frac{a}{q}\biggr)\int_{-R/qn}^{R/qn} M^{2}(\alpha)e(-m\alpha)
   d\alpha \\
   &+ \sum_{i,j} \frac{1}{\ph([q_{i},q])} \int_{-R/qn}^{R/qn}
   |M(\alpha)|d\alpha \cdot O\biggl(\frac{R}{q}(\log
   n)^{2}\Delta(n,[q_{j},q])\biggr) \\
   &+ O\biggl(\frac{R^{3}}{nq^{3}} (\log n)^{4}
   \Delta(n,[q_{2},q])\Delta(n,[q_{3},q])\biggr) \\
   &=: H_{a,q} (m) + \mathcal{O}_{1} + \mathcal{O}_{2},
\end{align*} 
say. Now
\begin{align*}
  \sum_{q_{2},q_{3}} \sum_{q\leq R} \sum_{\substack{ a<q \\ (a,q)=1
    }} \mathcal{O}_{1} &\ll \sum_{i,j} \sum_{q\leq R} \sum_{q_{i}}
  \frac{1}{\ph([q_{i},q])} \sum_{\substack{a<q \\ (a,q)=1}}
  \frac{R}{q} (\log n)^{2} \sum_{q_{j}} \Delta(n,[q_{j},q])\\
   &\ll \sum_{i,j} \sum_{q_{i}} \frac{\loglog n}{q_{i}} R (\log n)^{2} 
   \sum_{q_{j}} \sum_{q\leq R} \Delta(n,[q_{j},q]) \\
    &\ll (\log n)^{4} R \sum_{j} \sum_{h_{j}\leq RQ_{j}} \omega(h_{j})
    \Delta(n,h_{j})
\end{align*}
with
\begin{align*}
  \omega(h_{j}) &:= \sum_{q_{j}} \sum_{\substack{q\leq R \\
      [q_{j},q]=h_{j}}} 1 = \sum_{d_{j}\leq R} \sum_{q_{j}}
  \sum_{\substack{q\leq R \\ (q,q_{j})=d_{j} \\ qq_{j}=h_{j}d_{j} }} 1
  \\
  &\ll \sum_{d_{j}\leq R} \sum_{\substack{q\leq R \\ d_{j}|q}} 1 \ll
  R\log R \ll R\log n.
\end{align*}

So the $\mathcal{O}_{1}$-error term is
\begin{align*}
  &\ll R^{2}(\log n)^{5}\sum_{j} \sum_{h_{j}\leq RQ_{j}}
  \Delta(n,h_{j}) \ll R^{2} (\log n)^{5}\cdot \frac{n}{(\log
    n)^{\vartheta-B-6}} \\ &\ll n(\log n)^{3B-\vartheta+11} \ll n(\log
  n)^{-A-B-2} \ll n,
\end{align*}
again by using Bombieri-Vinogradov's Theorem and $\vartheta\geq A+4B+13$.

Now to $\mathcal{O}_{2}$. We have
\begin{align*}
  \sum_{q_{2},q_{3}} \sum_{q\leq R} \sum_{\substack{ a<q \\ (a,q)=1 }}
  \mathcal{O}_{2} &\ll  \sum_{q\leq R}\sum_{\substack{ a<q \\ (a,q)=1
    }} \frac{R^{3}}{nq^{3}}(\log n)^{4}\sum_{q_{2},q_{3}}
  \Delta(n,[q_{2},q]) \Delta(n,[q_{3},q]) \\
     &\ll  \frac{R^{3}}{n}(\log n)^{4} \sum_{\substack{ h_{2}\leq
         RQ_{2} \\ h_{3}\leq RQ_{3}}} \omega(h_{2},h_{3})
     \Delta(n,h_{2})\Delta(n,h_{3}) 
\end{align*}
with
\begin{align*}
  \omega(h_{2},h_{3}) &:= \sum_{q_{2},q_{3}} \sum_{\substack{q\leq R
      \\ [q_{i},q]=h_{i} \\ i=2,3}} \frac{1}{q^{2}} = \sum_{d_{2},d_{3}\leq R}
  \sum_{q_{2},q_{3}} \sum_{\substack{ q\leq R \\ (q_{i},q)=d_{i} \\
      q_{i}q = h_{i}d_{i} \\ i=2,3 }} \frac{1}{q^{2}} \\
    &\ll \sum_{d_{2},d_{3}\leq R} \sum_{\substack{q\leq R \\ [d_{2},d_{3}]|q }}
  \frac{1}{q^{2}} \ll \sum_{d_{2},d_{3}} \sum_{q\leq R}
  \frac{1}{q^{2}[d_{2},d_{3}]^{2}} \\ 
    &= \sum_{d_{2},d_{3}} \sum_{q\leq R}
    \frac{1}{q^{2}d_{2}^{2}d_{3}^{2}} (d_{2},d_{3})^{2} \ll
    \sum_{d_{3}\leq R} 1 \ll R,
\end{align*}
so the $\mathcal{O}_{2}$-error term is
\[
    \ll \frac{R^{4}}{n}(\log n)^{4} \Bigl( \sum_{h_{2}\leq RQ_{2}}
    \Delta(n,h_{2}) \Bigr) \Bigl( \sum_{h_{3}\leq RQ_{3}}
    \Delta(n,h_{3}) \Bigr) \ll n (\log n)^{6B-2\vartheta+16}\ll n,
\]
again by using Bombieri-Vinogradov's Theorem and $\vartheta\geq A+4B+13$.

Now there remains the main term.
Since
\[
    \int_{-R/qn}^{R/qn} M^{2}(\alpha) e(-m\alpha) d\alpha =
    m-1 + O\biggl(\frac{qn}{R}\biggr) \ll n
\]
for $q\leq R$ we can estimate it in the following way. It is
\begin{align*}
  H &:= \sum_{q_{2},q_{3}}\max_{\substack{m\leq n\\a_{2},a_{3}}}
    \sum_{q\leq R} \sum_{\substack{a<q \\ (a,q)=1}}
   \frac{(c_{2}c_{3})(a,q)}{\ph([q_{2},q])\ph([q_{3},q])}
   e\biggl(-m\frac{a}{q}\biggr) \int_{-R/qn}^{R/qn} M^{2}(\alpha)
   e(-m\alpha) d\alpha \\
   &\ll n \sum_{q_{2},q_{3}} \sum_{q\leq R} \frac{q(\log
     n)}{[q_{2},q][q_{3},q]} = n(\log n) \sum_{q_{2},q_{3}}
   \sum_{q\leq R} \frac{(q_{2},q)(q_{3},q)}{q_{2}qq_{3}} \\ 
   &\ll n (\log n) \sum_{a,b,c,d,e,f,g} \frac{dc\cdot db}{dace\cdot
     dabf \cdot dbcg} \ll n(\log n)^{7},
\end{align*}
where we substituted $q_{2}=dace$, $q_{3}=dabf$, $q=dbcg$ with
$a,b,c,d,e,f,g\leq n$, $d:=(q,q_{1},q_{3})$, and pairwise relatively
prime $a,b,c$ and $e,f,g$.

This shows the Proposition. \hfill$\square$

\subsubsection{Proof of Lemma \ref{lem3}}

Since the left hand side of Lemma 4 is
\begin{equation*}
  \ll \biggl( \sum_{q_{3}} q_{3} \max_{a_{3}} \sum_{q_{2}}q_{2} \max_{a_{2}}
  \sum_{m\leq n} |J_{2}^{\mathfrak{m}}(m)|^{2} \biggr)^{\!\hb},
\end{equation*}
it suffices to show that
\begin{equation*}
   \sum_{q_{3}} q_{3} \max_{a_{3}} \sum_{q_{2}}q_{2} \max_{a_{2}}
    \sum_{m\leq n} |J_{2}^{\mathfrak{m}}(m)|^{2} 
     \ll \frac{n^{3}}{(\log n)^{4A+32}}
\end{equation*}
for any $A>0$ in the required
regions for $Q_{2}$ and $Q_{3}$. The left hand side is
\begin{align*}
  &\sum_{q_{3}} q_{3} \max_{a_{3}} \sum_{q_{2}}q_{2} \max_{a_{2}}
    \sum_{m\leq n}
     \biggl|\int_{\mathfrak{m}} S_{2}(\alpha)S_{3}(\alpha)
     \:e(-m\alpha)d\alpha\biggr|^{2} \\
  &\leq \sum_{q_{3}} q_{3} \max_{a_{3}} \sum_{q_{2}}q_{2} \max_{a_{2}}
   \int_{\mathfrak{m}} |S_{2}(\alpha)S_{3}(\alpha)|^{2} d\alpha
\end{align*} 
by Bessel's inequality. Now
\begin{align*}
  |S_{2}(\alpha)|^{2} &= \sum_{\substack{m,m'\leq n \\ m\equiv
      m'\equiv a_{2}\: (q_{2})}} \Lambda(m) \Lambda(m')
  \:e(\alpha(m-m')) \\
   &= \sum_{\substack{|r|\leq n \\ r\equiv 0 \:(q_{2})}} e(\alpha r)
      \sum_{\substack{ m\leq n \\ m\equiv a_{2} \:(q_{2}) \\m-r\leq n
        }} \Lambda(m) \Lambda(m-r)   \\
   &=: \sum_{\substack{|r|\leq n \\ r\equiv 0 \:(q_{2})}} e(\alpha r)
   R(r;a_{2},q_{2}), 
\end{align*}
say, with $R(r;a_{2},q_{2})\ll \frac{n}{q_{2}}(\log n)^{2}$.

So the left hand side is
\begin{align*}
  &\ll n(\log n)^{2} \sum_{q_{3}\sim Q_{3}} q_{3} \max_{a_{3}}
   \sum_{q_{2}\sim Q_{2}}
   \sum_{\substack{|r|\leq n \\ r\equiv 0 \:(q_{2})}} 
   \biggl|\int_{\mathfrak{m}} |S_{3}(\alpha)|^{2}\:e(\alpha r)
    d\alpha\biggr| \\
  &\ll n(\log n)^{2} \sum_{q_{3}\sim Q_{3}} q_{3} \max_{a_{3}} 
    \sum_{0<|r|\leq n} \tau(|r|)
    \biggl| \int_{\mathfrak{m}} |S_{3}(\alpha)|^{2} 
    e(\alpha r)d\alpha\biggr| \\
  &+ n(\log n)^{2} \enspace Q_{2} \sum_{q_{3}\sim Q_{3}} q_{3} \max_{a_{3}} 
       \int_{0}^{1} |S_{3}(\alpha)|^{2} d\alpha.
\end{align*}
Now 
\[
   \int_{0}^{1}|S_{3}(\alpha)|^{2}d\alpha \ll \frac{n}{q_{3}}(\log n)^{2},
\]
so the second term is $\ll n^{2}(\log n)^{2}Q_{2}Q_{3}(\log n)^{2} \ll
n^{5/2}(\log n)^{4+\theta} \ll n^{3}(\log n)^{-A}$ and therefore in the
required bound.

The first term is
\begin{align*}
  &\ll n(\log n)^{2} \sum_{q_{3}} q_{3}\max_{a_{3}}
  \biggl(\sum_{0<|r|\leq n} \tau(|r|)^{2} \biggr)^{\!\h} 
  \biggl(\sum_{0<|r|\leq n} \biggl| \int_{\mathfrak{m}}
  |S_{3}(\alpha)|^{2} e(\alpha r) d\alpha \biggr|^{2} \biggr)^{\!\h} \\
  &\ll n^{3/2} (\log n)^{4} \sum_{q_{3}\sim Q_{3}} q_{3}\max_{a_{3}}
   \biggl( \int_{\mathfrak{m}} |S_{3}(\alpha)|^{4}
   d\alpha \biggr)^{\!\h} \\
  &\ll n^{3/2} (\log n)^{4} \biggl( \sum_{q_{3}\sim Q_{3}} 
    q_{3}^{2} \biggr)^{\!\h} \biggl( \sum_{q_{3}\sim Q_{3}}
    \max_{a_{3}} \int_{\mathfrak{m}} |S_{3}(\alpha)|^{4}
     d\alpha \biggr)^{\!\h} \\
  &\ll n^{3/2} (\log n)^{4} \biggl( \sum_{q_{3}\sim Q_{3}} q_{3}^{3}
    \max_{a_{3}} \int_{\mathfrak{m}} |S_{3}(\alpha)|^{4}
     d\alpha \biggr)^{\!\h}.
\end{align*}
Now here is the difficulty to show a nontrivial bound for the
expression in large brackets. It should be $\ll n^{3}/(\log n)^{C}$ 
for any large constant $C>0$ and large $Q_{3}$, but however one tries
to manage it, there is still some power of $Q_{3}$ left. We best can
give the bound
\begin{align*}
  \ll n^{3/2} (\log n)^{4} \biggl( \sum_{q_{3}\sim Q_{3}} q_{3}^{3}
    \max_{a_{3}} \max_{\alpha\in\mathfrak{m}} |S_{3}(\alpha)|^{2}
     \int_{0}^{1}|S_{3}(\alpha)|^{2}d\alpha \biggr)^{\!\h}.
\end{align*}

Now we need another Lemma to estimate $|S_{3}(\alpha)|^{2}$ for
$\alpha\in\mathfrak{m}$, it is the following.
\begin{lemma}
\label{lem4}
  For all $q_{3}\sim Q_{3}$, $(a_{3},q_{3})=1$ and
  $\alpha\in\mathfrak{m}$ we have $|S_{3}(\alpha)|^{2} \ll
  \frac{n^{2}}{q_{3}(\log n)^{C}}$ for $C=8A+2\theta+74$.
\end{lemma}

By using this we get for the above expression
\begin{align*}
  &\ll n^{3/2} (\log n)^{5} \biggl(\sum_{q_{3}\sim Q_{3}} Q_{3}
  \frac{n^{3}}{(\log n)^{C}} \biggr)^{\!\h} \\
  &\ll \frac{n^{3}}{(\log n)^{C/2-5}}Q_{3} \ll \frac{n^{3}}{(\log n)^{4A+32}}
\end{align*}
for $C= 8A+2\theta +74$ since $Q_{3}\leq (\log n)^{\theta}$.

So we see that $Q_{3}$ cannot be chosen as a power of $n$ using the
given method.

But this estimation shows Lemma \ref{lem3} for $Q_{2}\leq n^{\h}/(\log
n)^{\vartheta}$ and $Q_{3}\leq (\log n)^{\theta}$ as required.
\hfill$\square$

\zeile
\textbf{Proof of Lemma \ref{lem4}.} 

By Lemma 2 of A.\ Balog in \cite{c3} we have the validity of the following
assertion. For $C>0$ there exists a $D=D(C)>0$ such that for any
$\alpha\in\R$ with $\|\alpha-\frac{u}{v}\|<\frac{1}{v^{2}}$ with
integers $(u,v)=1$ and $(\log n)^{D}\leq v\leq \frac{n}{(\log n)^{D}}$
we have
\[
    \sum_{q_{3}\leq n^{1/3}/(\log n)^{D}} q_{3}
    \max_{(a_{3},q_{3})=1} |S_{3}(\alpha)|^{2} \ll \frac{n^{2}}{(\log n)^{C}},
\]
and since $Q_{3}\leq (\log n)^{\theta}\ll \frac{n^{1/3}}{(\log n)^{D}}$ also
\[
    \sum_{q_{3}\sim Q_{3}} q_{3}
    \max_{(a_{3},q_{3})=1} |S_{3}(\alpha)|^{2} \ll \frac{n^{2}}{(\log n)^{C}}.
\]

By Dirichlet's Approximation Theorem, for $\alpha\in\R$
and $B>0$ there exist integers $u,v$, $1\leq v\leq n/(\log n)^{B}$,
with $(u,v)=1$ and 
$\|\alpha-\frac{u}{v}\|<\frac{(\log n)^{B}}{vn}$,
and for $\alpha\in\mathfrak{m}$ it follows that $v\geq (\log n)^{B}$.

Therefore the conditions of Balog's Lemma are fulfilled if we take
$B\geq D(8A+2\theta+74)$, and it can be
applied then. It follows that for all $\alpha\in\mathfrak{m}$ we have
\[
   \sum_{q_{3}\sim Q_{3}} q_{3} \max_{a_{3}} |S_{3}(\alpha)|^{2}\ll
   \frac{n^{2}}{(\log n)^{8A+2\theta+74}},
\]
and so we have for all $q_{3}\sim Q_{3}$ and $(a_{3},q_{3})=1$ the
inequality
\[
   |S_{3}(\alpha)|^{2}\ll \frac{n^{2}}{q_{3}(\log n)^{8A+2\theta+74}},
\] 
since
\[
   |S_{3}(\alpha)|^{2} \ll \frac{1}{Q_{3}} \sum_{q_{3}\sim Q_{3}}
   q_{3} \max_{a_{3}} |S_{3}(\alpha)|^{2} \ll \frac{1}{Q_{3}}\cdot
   \frac{n^{2}}{(\log n)^{8A+2\theta+74}}.
\]
That shows Lemma \ref{lem4}. \hfill$\square$

\subsection{Proof of Theorem \ref{th2}}
\label{secsix}

Now we prove Theorem \ref{th2} in this last section.
Let $A,\theta, \vartheta>0$ be as in Theorem \ref{th2} and
let $n$ be odd and sufficiently large.

Besides $J_{3}(n)$ consider also 
\begin{equation*}
  R_{3}(n)=\sum_{\substack{p_{1},p_{2},p_{3} \\ p_{1}+p_{2}+p_{3}= n
    \\p_{i}\equiv a_{i} \;(q_{i}),\\ i=1,2,3}} \log p_{1} \log p_{2} \log p_{3}
\quad\text{ and }\quad
  r_{3}(n)=\sum_{\substack{p_{1},p_{2},p_{3} \\ p_{1}+p_{2}+p_{3}= n
    \\p_{i}\equiv a_{i} \;(q_{i}),\\i=1,2,3}} 1.
\end{equation*}
Then we have
\begin{equation*}
   |R_{3}(n)-J_{3}(n)|\leq (\log n)^{3}W,  
\end{equation*}
where $W$ denotes the number of solutions of $p^{l}+q^{j}+r^{k}=n$,
with $p,q,r$ prime and where $l$, $j$ or $k$ are at least $2$,
and $p^{l}\equiv a_{1} \:(q_{1})$, $q^{j}\equiv a_{2} \:(q_{2})$,
$r^{k}\equiv a_{3} \:(q_{3})$. Now four cases occur: For $i=1,2,3,4$
let $W_{(i)}$ be the number of solutions in case $(i)$, namely
$(1)\: l,j\geq 2$, $(2)\: l=1,j\geq 2$, $(3)\: l\geq 2,j=1$,
$(4)\: l=j=1, k\geq 2$.

In case $(1)$ there are at most $O(\sqrt{n})$ many possibilities 
for $p^{l},q^{j}\leq n$, so $W_{(1)}\ll n$ and we have $\sum_{q_{1},q_{2},q_{3}}
\max_{a_{1},a_{2},a_{3}} W_{(1)}\ll \frac{n^{2}}{(\log n)^{2\vartheta-\theta}}
\ll\frac{n^{2}}{(\log n)^{A+3}}$ since $\vartheta>\theta+A+3$.

In case $(4)$ we have at most $O(\sqrt{n})$ many possibilities for
$r^{k}\leq n$ and $\ll \frac{n}{q_{2}}$ many for $q$,
so $W_{(4)}\ll \frac{n^{3/2}}{q_{2}}$ and we get $\sum_{q_{1},q_{2},q_{3}}
\max_{a_{1},a_{2},a_{3}} W_{(4)}\ll Q_{1}Q_{3}n^{3/2}
\ll\frac{n^{2}}{(\log n)^{\vartheta-\theta}}
\ll \frac{n^{2}}{(\log n)^{A+3}}$ since $\vartheta>\theta+A+3$.

The same estimation comes of course analogously with $W_{(2)}$ in case
$(2)$.

In case $(3)$ we consider the number
\begin{equation*}
  \#\{p^{l}\leq n;\: l\geq 2, p^{l}\equiv a_{1} \:(q_{1}) \} \leq
  \sum_{\substack{m\leq n\\m\equiv a_{1}(q_{1})}}
  \Lambda(m)(1-\mu^{2}(m))=: N(a_{1},q_{1})
\end{equation*}
in the context of section \ref{secthree}, with $b_{m}:=
\Lambda(m)(1-\mu^{2}(m))$. Then $W_{(3)}\ll N(a_{1},q_{1})\cdot
\frac{n}{q_{2}}$, and by application of Lemma \ref{lem1} we get
\begin{align*}
  \sum_{q_{1},q_{2},q_{3}} \max_{a_{1},a_{2},a_{3}} W_{(3)}
   \ll nQ_{3}\sum_{q_{1}}  \max_{a_{1}} N(a_{1},q_{1})
  &\ll nQ_{3}\Bigl(\sum_{q_{1}} q_{1} \max_{a_{1}} N(a_{1},q_{1})^{2} \Bigr)^{\h} \\
  &\ll nQ_{3}\Bigl(\frac{n^{2}}{H} +n^{3/2}H \Bigr)^{\h} \log n
\end{align*}
since
\begin{equation*}
  \sum_{m\leq n} |b_{m}|^{2}=\sum_{\substack{m\leq n\\m\equiv a_{1} \,(q_{1})}}
  \Lambda^{2}(m) (1-\mu^{2}(m))^{2} \ll\sum_{\substack{p^{k}\leq n \\
      k\geq 2}} (\log p)^{2} \ll \sqrt{n} \log n
\end{equation*}
and $Q_{1}\leq \sqrt{n}$.

If we choose the parameter $H$ as $H:= \frac{n^{\h}}{(\log n)^{2A+6}Q_{3}^{2}}$
we get further
\begin{align*}
  \sum_{q_{1},q_{2},q_{3}} \max_{a_{1},a_{2},a_{3}} W_{(3)}
   &\ll nQ_{3}\Bigl( n^{3/2} (\log n)^{2A+6}Q_{3}^{2}
  + \frac{n^{2}}{Q_{3}^{2}(\log n)^{2A+6}} \Bigr)^{\h} \\
  &\ll n\cdot n^{3/4}Q_{3}^{2}(\log n)^{A+3} + \frac{n^{2}}{(\log n)^{A+3}}
   \ll \frac{n^{2}}{(\log n)^{A+3}}.
\end{align*}

So we get
\begin{equation*}
  \sum_{q_{1},q_{2},q_{3}} \max_{a_{1},a_{2},a_{3}} W
     \ll \frac{n^{2}}{(\log n)^{A+3}}.
\end{equation*}

Therefore it follows from Theorem \ref{Th1}:
\begin{align*}
  &\sum_{q_{3}} \max_{a_{3}} 
   \sum_{q_{2}} \max_{a_{2}} \sum_{q_{1}} \max_{a_{1}} \Bigl|R_{3}(n) -
  \frac{n^{2}\mathcal{S}_{3}(n)}{2\ph(q_{1})\ph(q_{2})\ph(q_{3})} \Bigr| \\
   \leq & \sum_{q_{3}} \max_{a_{3}} \sum_{q_{2}} \max_{a_{2}} 
    \sum_{q_{1}} \max_{a_{1}} \Bigl|R_{3}(n) - J_{3}(n)\Bigr| \\
      &+ \sum_{q_{3}} \max_{a_{3}} \sum_{q_{2}} \max_{a_{2}} 
     \sum_{q_{1}} \max_{a_{1}} 
  \Bigl|J_{3}(n) -  \frac{n^{2}\mathcal{S}_{3}(n)}%
                         {2\ph(q_{1})\ph(q_{2})\ph(q_{3})}\Bigr| \\
  \ll &\sum_{q_{1},q_{2},q_{3}} \max_{a_{1},a_{2},a_{3}} W(\log n)^{3} 
   + \frac{n^{2}}{(\log n)^{A}}  \ll  \frac{n^{2}}{(\log n)^{A}}.
\end{align*}
So the formula of Theorem \ref{Th1} holds also for $R_{3}(n)$ instead of $J_{3}(n)$.

Now let $q_{3}\leq Q_{3}= (\log n)^{\theta}$ and $(a_{3},q_{3})=1$ be
fixed.

For given $q_{2}$ and admissible $a_{2}$ consider
\begin{equation*}
  \mathcal{Q}_{1}:=\{ q_{1}\leq Q_{1};\: \exists\: a_{1}\text{ adm.}:
  R_{3}(n)=0 \},\quad E_{1}:=\#\mathcal{Q}_{1},
\end{equation*}
and
\begin{equation*}
  \mathcal{Q}_{2}:=\{ q_{2}\leq Q_{2};\: \exists\: a_{2}\text{ adm.}:
     E_{1}\geq Q_{1}(\log n)^{-A} \},\quad E_{2}:=\#\mathcal{Q}_{2}.
\end{equation*}
We have $\mathcal{S}_{3}(n)\gg 1$ if it is positive (see the formula
for it as Euler product), so we have
\begin{align*}
E_{2}\cdot \frac{Q_{1}}{(\log n)^{A}}&\cdot
\frac{n^{2}}{Q_{1}Q_{2}Q_{3}} \\
&\leq \sum_{q_{2}\in\mathcal{Q}_{2}}\; \max_{\substack{a_{2}\text{
      adm.}\\E_{1}\geq\frac{Q_{1}}{(\log n)^{A}}}} \;
\sum_{q_{1}\in\mathcal{Q}_{1}}\; \max_{\substack{a_{1}\text{ adm.}\\
  R_{3}(n)=0}}\;
\biggl|\frac{n^{2}\mathcal{S}_{3}(n)}{2\ph(q_{1})\ph(q_{2})\ph(q_{3})}
 \biggr| \\ &\ll \frac{n^{2}}{(\log n)^{2A+\theta}}
\end{align*}
by Theorem \ref{Th1}, and it follows that $E_{2}\ll Q_{2}(\log n)^{-A}$.

So for almost all $q_{2}$ and all admissible $a_{2}$ we have
that $E_{1}< Q_{1}(\log n)^{-A}$, that means that for almost all $q_{1}$
and all admissible $a_{1}$ it holds that $R_{3}(n)>0$.
Since $r_{3}(n)\geq \frac{R_{3}(n)}{(\log n)^{3}}$, it follows that 
$r_{3}(n)$ is positive, too, so Theorem \ref{th2} follows.
\hfill$\square$

%
%
%


\end{document}